\documentclass[final,3p,times]{elsarticle}

\makeatletter
\def\ps@pprintTitle{%
 \let\@oddhead\@empty
 \let\@evenhead\@empty
 \def\@oddfoot{}%
 \let\@evenfoot\@oddfoot}
\makeatother

\usepackage{amssymb}
\usepackage{amsthm}

\makeatletter
\newcommand*{\rom}[1]{\expandafter\@slowromancap\romannumeral #1@}
\makeatother

\newcommand{\calP}{\ensuremath{\mathcal{P}}}
\newcommand{\calE}{\ensuremath{\mathcal{E}}}
\newcommand{\calN}{\ensuremath{\mathcal{N}}}
\newcommand{\calM}{\ensuremath{\mathcal{M}}}
\newcommand{\calC}{\ensuremath{\mathcal{C}}}
\newcommand{\calT}{\ensuremath{\mathcal{T}}}
\newcommand{\calI}{\ensuremath{\mathcal{I}}}
\newcommand{\calD}{\ensuremath{\mathcal{D}}}
\newcommand{\calX}{\ensuremath{\mathcal{X}}}
\newcommand{\calF}{\ensuremath{\mathcal{F}}}
\newcommand{\calL}{\ensuremath{\mathcal{L}}}
\newcommand{\calA}{\ensuremath{\mathcal{A}}}
\newcommand{\calU}{\ensuremath{\mathcal{U}}}
\newcommand{\calZ}{\ensuremath{\mathcal{Z}}}
\newcommand{\calB}{\ensuremath{\mathcal{B}}}
\newcommand{\calQ}{\ensuremath{\mathcal{Q}}}
\newcommand{\calR}{\ensuremath{\mathcal{R}}}
\newcommand{\calS}{\ensuremath{\mathcal{S}}}
\newcommand{\calG}{\ensuremath{\mathcal{G}}}

\newcommand{\Mtens}{\underline{\calM_h}}
\newcommand{\Xtens}{\underline{\calX}}
\newcommand{\Atens}{\underline{\calA_h}}
\newcommand{\Ahattens}{\underline{\hat\calA_h}}
\newcommand{\Ctens}{\underline{\calC_h}}
\newcommand{\Ttens}{\underline{\calT_h}}
\newcommand{\Ztens}{\underline{\calZ}}
\newcommand{\Btens}{\underline{\underline{\calB}}}
\newcommand{\Dtens}{\underline{\calD}}
\newcommand{\Ftens}{\underline{\calF}}
\newcommand{\Fhattens}{\underline{\hat\calF}}
\newcommand{\Utens}{\underline{\calU}}
\newcommand{\Qtens}{\underline{\calQ_h}}

\newcommand{\XtensP}{\underline{\calX_h^P}}

\newcommand{\MRB}{\underline{\calM_{\rm RB}}}
\newcommand{\ARB}{\underline{\calA_{\rm RB}}}
\newcommand{\CRB}{\underline{\calC_{\rm RB}}}
\newcommand{\TRB}{\underline{\calT_{\rm RB}}}
\newcommand{\XRB}{\underline{Z_{\rm RB}}}
\newcommand{\QRB}{\underline{\calQ_{\rm RB}}}

\newcommand{\Bvect}{\underline{\calB}}

\newcommand{\RR}{\ensuremath{\mathbb{R}}}
\newcommand{\PP}{\ensuremath{\mathbb{P}}}
\newcommand{\II}{\ensuremath{\mathbb{I}}}
\newcommand{\UU}{\ensuremath{\mathbb{U}}}
\newcommand{\ZZ}{\ensuremath{\mathbb{Z}}}
\newcommand{\CC}{\ensuremath{\mathbb{C}}}

\newcommand{\bM}{\ensuremath{\mathbf{M}}}
\newcommand{\bN}{\ensuremath{\mathbf{N}}}
\newcommand{\bD}{\ensuremath{\mathbf{D}}}

\def\ii{\textrm{i}}
\def\ntt{n^\textrm{t-t}}
\def\ufreq{\hat{u}}
\def\mufreq{\tilde{\mu}}
\def\Pfreq{\tilde{\calP}}
\def\ufe{u_h}
\def\ufef{\hat{u}_h}
\def\ufefd{u_{h,\Delta t}}
\def\dufefd{\dot{u}_{h,\Delta t}}
\def\ddufefd{\ddot{u}_{h,\Delta t}}
\def\upr{\hat\UU_{h,\bD}}
\def\upri{\hat U_{h,\bD,i}}
\def\upro{\hat U_{h,\bD,1}}
\def\uprn{\hat U_{h,\bD,N}}
\def\urb{U_{h,\bD,N}}
\def\urbfd{U_{h,\bD,N,\Delta t}}
\def\urbfderr{U_{h,\bD,N,2\Delta t}}
\def\durbfd{\dot{U}_{h,\bD,N,\Delta t}}
\def\ddurbfd{\ddot{U}_{h,\bD,N,\Delta t}}

\def\afreq{\hat{a}}
\def\ffreq{\hat{f}}

\def\murefC{\mu^{{\rm ref},C}}
\def\mureffC{\hat\mu^{{\rm ref},C}}
\def\mureffP{\hat\mu^{{\rm ref},P}}

\newcommand{\AtensPmu}{\underline{\calA_{h\; \mufreq^P}^P}}
\newcommand{\BtensPmu}{\underline{\calB_{h\; \mufreq^P}^P}}
\newcommand{\DtensPmu}{\underline{\calD_{h\; \mufreq^P}^P}}
\newcommand{\QtensPmu}{\underline{\calQ_{h\; \mufreq^P}^P}}
\newcommand{\PtensPmu}{\underline{\calP_{h\; \mufreq^P}^P}}
\newcommand{\Ptenstemp}{\underline{\calP_{h\; {\rm temp}}^P}}

\newcommand{\GtensR}{\underline{\calG_{h\; R}^P}}
\newcommand{\GtensS}{\underline{\calG_{h\; S}^P}}

\graphicspath{{fig_PR_RBC/}}
\usepackage{amsmath,amsfonts,stmaryrd,latexsym,amstext}
\usepackage{bbm}
\usepackage{float}
\usepackage{algorithm}
\usepackage[noend]{algpseudocode}
\usepackage{xcolor}

\journal{Computer Methods in Applied Mechanics and Engineering}

\begin{document}


\begin{frontmatter}



 \title{A Two-Level Parameterized Model-Order Reduction Approach for Time-Domain Elastodynamics \tnoteref{t1}}
 \tnotetext[t1]{This work was supported by the Office of Naval Research [N00014-17-1-2077]; and the Army Research Office [W911NF1910098]}


\author[1]{Mohamed Aziz BHOURI\corref{cor1}\fnref{fn1}%
  }
\ead{mohamedazizbhouri@gmail.com}
\author[1]{Anthony T. PATERA\corref{cor2}
}
\ead{patera@mit.edu}

 \cortext[cor1]{Corresponding author}
 \fntext[fn1]{Present address: 3401 Walnut St, Wing A, Office 536, Philadelphia, PA 19104, USA  \\ \\ \\ Preprint submitted to Computer Methods in Applied Mechanics and Engineering}
 \address[1]{Department of Mechanical Engineering, Massachusetts Institute of Technology, 77 Massachusetts Avenue, Cambridge, MA 02139 USA}

\begin{abstract}
We present a two-level parameterized Model Order Reduction (pMOR) technique for the linear hyperbolic Partial Differential Equation (PDE) of time-domain elastodynamics. In order to approximate the frequency-domain PDE, we take advantage of the Port-Reduced Reduced-Basis Component (PR-RBC) method to develop (in the offline stage) reduced bases for subdomains; the latter are then assembled (in the online stage) to form the global domains of interest. The PR-RBC approach reduces the effective dimensionality of the parameter space and also provides flexibility in topology and geometry. In the online stage, for each query, we consider a given parameter value and associated global domain. In the first level of reduction, the PR-RBC reduced bases are used to approximate the frequency-domain solution at selected frequencies. In the second level of reduction, these instantiated PR-RBC approximations are used as surrogate truth solutions in a Strong Greedy approach to identify a reduced basis space; the PDE of time-domain elastodynamics is then projected on this reduced space. We provide a numerical example to demonstrate the computational capability and assess the performance of the proposed two-level approach.
\end{abstract}



\begin{keyword}
Model order reduction \sep domain decomposition \sep parametrized partial differential equations \sep elastodynamics 


\end{keyword}

\end{frontmatter}



\section{Introduction}
\label{sec:intro}

Model Order Reduction (MOR) for time-dependent problems has received a great deal of attention in the reduced-order modeling community. These methods are particularly appropriate for many-query applications. For example, MOR methods for the linear time-domain elastodynamics equation are of great interest in the context of Simulation Based Classification (SBC) for Structural Health Monitoring (SHM) based on time-domain cross-correlation functions \cite{Huo2016,Yang2009}; the latter requires the construction of large (synthetic) training data sets of the time-domain response of mechanical structures under ambient localized excitations.

Existing MOR techniques for time-dependent problems include Proper Orthogonal Decomposition (POD) approaches \cite{KunischVolkwein,KunischVolkwein2,RathinamPetzold}, Greedy methods \cite{Grepletal,GreplPatera}, hybrid approaches combining POD (in time) and Greedy procedures (in parameter space) \cite{Grepl}, and space-time approaches \cite{MasaThesis}. A second set of approaches are based on the frequency-time duality: the reduced space is constructed in the frequency domain, and the time-domain equation is then projected on the reduced space. One such technique, used extensively within the control and dynamics community, is interpolatory model reduction \cite{Antoulasetal,Beattieetal}. A similar reduced basis method for the time-domain heat equation and wave equation has also been developed \cite{HuynhetalCR}. However, most existing methods can not address many parameters, local excitation, or topology variation, all of which are important in (say) the SBC-SHM context.

Recent advances in MOR offer new opportunities for the development of more efficient approaches for linear time-domain PDEs. In particular, domain decomposition (DD)-based MOR techniques can reduce the effective dimensionality of the parameter space considered in the construction of the local reduced bases. Among those techniques, the Port-Reduced, Reduced-Basis Component (PR-RBC) method \cite{HuynhetalESAIM,HuynhetalCMAME,EftangPatera,Smetana,SmetanaPatera} provides both rapid response and also flexibility in topology and geometry.  The PR-RBC is essentially a combination of the Component Model Synthesis (CMS) technique \cite{CraigBampton,Hurty,HetmaniukLehoucq} | as regards components and ports | and the Reduced Basis method \cite{NoorPeters,Almrothetal,Rozzaetal,Barraultetal} | as regards bubbles and in particular parametric treatment. The first synthesis of CMS and RB is the Reduced Basis Element method (RBE) \cite{MadayRonquist}; PR-RBC may be viewed as a Reduced Basis Element method for a particular (Static Condensation \cite{Wilson}) choice for the interface treatment and particular strategies for port mode training \cite{EftangPatera} and bubble-mode training \cite{Veroyetal}.

In the PR-RBC approach, we define a library of parameterized archetype components: each of these components is characterized by ``local" parameters, a reference finite element (FE) mesh, and local ports for interconnection; the latter induce an associated library of reference ports. For every archetype component, reduced bases are built to approximate the solution inside the domain due to any sources; the latter are referred to as the ``Reduced Bubble Spaces for Inhomogeneity". For every reference port, a low-dimensional space | a set of ``port modes" | is built to approximate the behavior of the solution over the reference port by a component-pairwise training procedure \cite{EftangPatera,SmetanaPatera}. For each such port mode, a low-dimensional space, denoted Reduced Bubble Space for Port Mode Liftings, is constructed to approximate the port mode lifting. The space construction proceeds in an offline stage. Then, in the online stage, instantiated components are assembled into a global system and associated RB spaces are considered; the global parameter is prescribed and the linear equations are solved by static condensation.

Our goal here is to develop an efficient RB method for long-time integration of the linear time-domain elastodynamics PDE for large geometric domains with localized excitations and relatively many parameters. Towards that end, we develop a two-level parameterized model-order reduction approach, henceforth referred to as two-level PR-RBC method, by extending frequency-to-time-domain methods to incorporate component approaches \cite{HuynhetalESAIM,HuynhetalCMAME,EftangPatera,Smetana,SmetanaPatera}. The standard PR-RBC method is applied to the frequency-domain equation, where an augmented parameter set is considered that also includes the frequency. In an offline stage, the necessary spaces are formed for a given library of components. Then a global structure is formed as an assembly of instantiated archetype components. The PR-RBC approximation is then used to compute the global frequency-domain solutions at well-selected frequencies. This corresponds to the first-level reduction. Thereafter, these solutions are taken as high-fidelity approximations | a surrogate for the FE ``truth" | to form a (final) reduced basis to approximate the global time-domain solution. This corresponds to the second level of reduction. In this work, we choose a strong greedy procedure to form the (final) reduced basis from the PR-RBC snapshots: this approach is justified by the efficiency of the PR-RBC method. Our two-level approach arises from the incorporation of the PR-RBC procedure | applied to the frequency-domain PDE | into the approximation of the time-domain PDE. Therefore, the two levels of reduction are completely distinct: the first level consists of the PR-RBC approximation of the frequency-domain solution, while the second level consists of RB approximation of the time-domain solution. In contrast, earlier two-level model reduction approaches rely on standard RB methods for both levels \cite{Eftangetal}.

This paper is organized as follows. In Section \ref{sec:elast_form}, we introduce the weak form of the time-domain elastodynamics equation of interest and the corresponding frequency-domain equation. In Section \ref{sec:FEapp}, we present the finite-element and finite-difference discretizations which constitute our ``truth" approximation | the point of departure for subsequent model order reduction. An overview of the PR-RBC method is given in Section \ref{sec:PRRBC}. The proposed two-level PR-RBC method is developed in Section \ref{sec:2RB}: we present the formulation and provide an operation count. Finally, in order to demonstrate the capability and assess the performance of the proposed technique, a numerical example is presented in Section \ref{sec:num_ex}.

\section{Elastodynamics Formulation}
\label{sec:elast_form}

The two-level PR-RBC method proposed in this work can be applied to any linear time-domain PDE which admits an affine representation of the parameter. The latter can be recovered by means of EQP \cite{PateraYano} or EIM \cite{Barraultetal,Grepletal} if needed. For sake of clarity, we restrict ourselves to the PDE of linear elastodynamics in this work.

\subsection{Time-Domain Equation}

Let $\Omega\subset\RR^d$, $d=2$ be a bounded domain, with boundary $\partial\Omega$. The proposed method can be naturally applied to the 3D-elastodynamics case, but for purposes of presentation we limit ourselves to the 2D case in this work. The boundary $\partial\Omega$ is assumed to be partitioned into $\Gamma^D$ and $\Gamma^N$, such that Dirichlet boundary conditions are imposed on $\Gamma^D$, while natural boundary conditions are satisfied on $\Gamma^N$. Without loss of generality, the Dirichlet boundary conditions are assumed to be homogeneous, and the non-essential boundary conditions are assumed to be of Neumann type; more details on treating non-homogeneous Dirichlet boundary conditions can be found in \cite{BhouriThesis}. The simulation time interval is noted $[0,T_{\rm final}]$, $T_{\rm final}>0$. Let ${X}\equiv\{v\in [H^1({\Omega})]^{d} \ | \ v|_{{\Gamma}^{D}}=0\}$ be the Hilbert space of admissible real-valued functions; ${X}$ is imbued with inner product $(w,v)_{{X}} \equiv \int_{{\Omega}} \nabla w \cdot \nabla v + w v \, dV$ and induced norm $\| w \|_{{X}} \equiv \sqrt{(w,w)_{{X}}}$. The problem parameterization is denoted ${\mu}\in{\calP}$, where ${\calP}\in\RR^{{n}_P}$ is a suitable compact set.


Let $m(\cdot,\cdot;\cdot) ; c(\cdot,\cdot;\cdot) ; a(\cdot,\cdot;\cdot) :\big(H^1({\Omega})\big)^d\times \big(H^1({\Omega})\big)^d\times\calP\rightarrow\mathbb{R}$ be the bilinear forms corresponding to the mass, damping and stiffness terms of the 2D-elastodynamics equation respectively. These bilinear forms are defined as follows:
\begin{equation}
m(w,v;\mu)\equiv\rho\;\int_{\Omega}w\cdot  v \ dx \ ,
\end{equation}
\begin{equation}
a(w,v;\mu)\equiv\frac{\nu\; E}{(1+\nu)\;(1-2\nu)}\;\int_{\Omega}\frac{\partial w_i}{\partial x_j}\frac{\partial {v_k}}{\partial x_l}\delta_{ik}\delta_{jl} \ dx + \frac{E}{2(1+\nu)}\;\int_{\Omega}\frac{\partial w_i}{\partial x_j}\frac{\partial {v_k}}{\partial x_l}(\delta_{ik}\delta_{jl}+\delta_{il}\delta_{jk}) \ dx \ , \label{eq:a_bil_form_time_domain}
\end{equation}
\begin{equation}
c(w,v;\mu)\equiv\alpha_{\rm Ray}\; m(w,v;\mu) + \beta_{\rm Ray}\; a(w,v;\mu) \ , \label{eq:c_bil_form_time_domain}
\end{equation}
\noindent where in equation (\ref{eq:a_bil_form_time_domain}) we use the convention of summation over repeated indices, $\nu$ denotes the Poisson ratio, $\rho>0$ is the material density, $\alpha_{\rm Ray}>0$ and $\beta_{\rm Ray}>0$ are the Rayleigh damping coefficients, and $E>0$ is the Young's modulus. Here $f(\cdot,\cdot;\cdot) :\big(H^1({\Omega})\big)^d\times[0,T_{\rm final}]\times\calP\rightarrow\mathbb{R}$ refers to linear form corresponding to the Neumann boundary conditions imposed on $\Gamma^N$. In order to realize an efficient offline-online decomposition, the bilinear and linear forms are assumed to have an affine dependence on the parameter $\mu$. We further assume that the linear form $f(\cdot,\cdot;\cdot)$ satisfies the following space-time separation of variable:
\begin{equation}
\label{eq:fxft}
f(v,t;\mu) = f_t(t;\mu_t)\; f_x(v;{\mu_0}) \ , \forall v\in X \ , \forall t\in[0,T_{\rm final}] \ , \forall\mu\in\calP \ ,
\end{equation}
\noindent where $\mu_t$ denotes the parameters governing $f_t(t,\cdot)$, and ${\mu_0}$ refers to all remaining parameters, with $\mu=(\mu_t,\mu_0)$.

The variational formulation of the elastodynamics equation then reads as follows: Find $u(t\in[0,T_{\rm final}];\mu)$ such that $\forall t\in[0,T_{\rm final}]$ ,
\begin{equation}
m\Big(\frac{\partial^2 u(t;\mu)}{\partial t^2},v;\mu\Big)+c\Big(\frac{\partial u(t;\mu)}{\partial t},v;\mu\Big)+a\Big(u(t;\mu),v;\mu\Big)=f(v,t;\mu) \ , \forall v\in X \ , \forall t\in[0,T_{\rm final}]\ , \label{eq:VarFormGlobElast}
\end{equation}
\begin{equation}
u(t;\mu)=0 \ , {\rm on} \ \Gamma^D \ , \forall t\in[0,T_{\rm final}] \ ,
\label{eq:BCDirichletGlob}
\end{equation}
\noindent and
\begin{equation}
u(t=0;\mu)=0 \ ; \frac{\partial u}{\partial t}\Big(t=0;\mu\Big)=0 \ . 
\end{equation}
\noindent Extension to non-zero initial conditions can also be considered.

\subsection{Frequency-Domain Equation}
\label{subsec:FreqEq}

Let $\hat X\equiv\{v \ | \ v = w + \ii y \ , w,y\in X\}$ be the Hilbert space of admissible complex-valued functions; ${\hat X}$ is imbued with inner product $(w,v)_{\hat{X}} \equiv \int_{{\Omega}} \nabla w \cdot \overline{\nabla v} + w \bar v \, dV$ and induced norm $\| w \|_{\hat{X}} \equiv \sqrt{(w,w)_{\hat{X}}}$, where $\overline{\ \cdot \ }$ refers to the complex conjugate operator. The complex problem parameterization is denoted $\mufreq=(\mu,\omega)\in\Pfreq$, which corresponds to the real-valued problem parameter concatenated with the angular frequency $\omega$ as an additional parameter; here $\Pfreq\in\RR^{{n}_P+1}$ refers to the augmented compact parameter set.


We assume that we have non-zero damping such that $\int\limits_0^{\infty}|u(t;\mu)|dt<\infty$, and write $u(t;\mu)=\Re \{\ufreq\big[ \mufreq=(\mu,\omega) \big] e^{ \ii \omega t} \}$, where $\Re$ refers to real part. It follows that $\ufreq(\mufreq)$ satisfies the variational formulation of the Helmholtz equation: Find $\ufreq(\mufreq)\in\hat X$ such that:
\begin{equation}
\afreq(\ufreq(\mufreq),v;\mufreq)=\ffreq(v;\mufreq) \ , \forall v\in\hat X \ , \label{eq:VarFormGlobHom}
\end{equation}
\noindent where $\afreq(\cdot,\cdot;\cdot):\hat X\times \hat X\times\Pfreq\rightarrow\mathbb{C}$ is a sesquilinear form given by
\begin{equation}
\afreq(\cdot,\cdot;\cdot)=-\omega^2\; \hat m(\cdot,\cdot;\cdot)+\ii\omega\; \hat c(\cdot,\cdot;\cdot)+ \hat k(\cdot,\cdot;\cdot) \ .
\end{equation}
\noindent Here $\hat m(\cdot,\cdot;\cdot)$, $\hat c(\cdot,\cdot;\cdot)$ and $\hat k(\cdot,\cdot;\cdot)$ are defined as
\begin{equation}
\hat m(w,v;\mufreq)=\rho\;\int_{\Omega}w\cdot \bar v \ dx \ ,
\end{equation}
\begin{equation}
\hat k(w,v;\mufreq)=\frac{\nu\; E}{(1+\nu)\;(1-2\nu)}\;\int_{\Omega}\frac{\partial w_i}{\partial x_j}\frac{\partial \overline{v_k}}{\partial x_l}\delta_{ik}\delta_{jl} \ dx + \frac{E}{2(1+\nu)}\;\int_{\Omega}\frac{\partial w_i}{\partial x_j}\frac{\partial \overline{v_k}}{\partial x_l}(\delta_{ik}\delta_{jl}+\delta_{il}\delta_{jk}) \ dx \ , \label{eq:k_bil_form_freq_domain}
\end{equation}
\begin{equation}
\hat c(w,v;\mufreq)=\alpha_{\rm Ray}\; \hat m(w,v;\mufreq) + \beta_{\rm Ray}\; \hat a(w,v;\mufreq) \ ,
\end{equation}
\noindent where in equation (\ref{eq:k_bil_form_freq_domain}) we use the convention of summation over repeated indices and $ \ffreq(\cdot;\cdot):\hat X\times\Pfreq\rightarrow\mathbb{R}$ is a continuous anti-linear form associated to the Laplace transform of the $f(\cdot,\cdot;\cdot)$ term. 

Since $m(\cdot,\cdot;\mu)$, $c(\cdot,\cdot;\mu)$, $a(\cdot,\cdot;\mu)$ and $f(\cdot,\cdot;\mu)$ are assumed to have an affine dependence on the parameter $\mu$, it follows that $\afreq(\cdot,\cdot;\mufreq)$ and $\ffreq(\cdot;\mufreq)$ also have an affine dependence on the parameter $\mufreq$. Moreover, since $f(\cdot,\cdot;\cdot)$ is assumed to satisfy the space-time separation of variable (\ref{eq:fxft}), and since $f_t(t;\mu_t)\rightarrow u(t;\mu)$ is a linear time invariant system, we can write:
\begin{equation}
u\Big(t;\mu=(\mu_t,\mu_0)\Big)=\int_{-\infty}^{\infty} f_t(\tau;{\mu_t})\; h(t-\tau;\mu_0)\; d\tau \ ,
\end{equation}
\begin{equation}
\ufreq\Big(\mufreq=(\omega,\mu_t,\mu_0)\Big)=\hat h(\omega,\mu_0)\;\hat f_t(\omega,\mu_t) \ ,
\end{equation}
\noindent where $\hat f_t(\omega,\mu_t)=\calL\calT[f_t(\cdot;\mu_t)]$ ($\calL\calT$ refers to the Laplace transform operator) and $h(\cdot;\mu_0)$, and $\hat h(\cdot,\mu_0)$ denote the space-dependent time-domain and frequency-domain representations of the transfer function respectively. If $\hat f_t(\omega,\mu_t)=1$, then $\ufreq\Big((\omega,\mu_t,\mu_0)\Big)=\hat h(\omega,\mu_0)$. Therefore, if we take $\hat f_t(\omega,\mu_t)=1$ and perform a model order reduction approach, the reduced basis will learn $\hat h(\omega,\mu_0)$ and thus we can well approximate $\ufreq\Big((\omega,\mu_t,\mu_0)\Big)=\hat h(\omega,\mu_0)\;\hat f_t(\omega,\mu_t)$, for any given $\hat f_t(\omega;\mu_t)$. Since $u(t;\mu)=\calL\calT^{-1}[\ufreq(\mufreq)]$, we can also well approximate $u(\cdot;\mu)$. For the remainder of this work, any frequency-domain problem will be considered to have an anti-linear form which is constant and equal to $1$, independent of the time-dependence of the linear form of the time-domain problem.

\section{Finite Element Approximation}
\label{sec:FEapp}

\subsection{Finite Element Discretization}
In order to approximate equations (\ref{eq:VarFormGlobElast}) and (\ref{eq:VarFormGlobHom}), we consider a suitably refined finite element (FE) Galerkin approximation: a triangulation ${\calT}^h$ for domain $\Omega$; associated conforming FE approximation spaces $X_h^0\subset \big(H^1(\Omega)\big)^d$ for real-valued functions and $\hat X_h^0$ for complex-valued functions, both of dimension $\calN^0_h$. Let $X_h\equiv X\cap X_h^0$ and $\hat X_h\equiv \hat X\cap \hat X_h^0$, and let $\{\varphi_j\}_{j=1,\ldots,\calN_h}$ denote the associated (real) standard FE nodal basis; here $\calN_{h}$ is the dimension of $X_h$ and $\hat X_h$.


\subsubsection{Time-Domain equation}
\label{sssec:time}

The FE approximation $\ufe(t\in[0,T_{\rm final}];\mu)$ to $u(t\in[0,T_{\rm final}];\mu)$ can be obtained by projecting equation (\ref{eq:VarFormGlobElast}) on $X_h$: $\ufe(t;\mu)\in X_{h}\ , \forall t\in[0,T_{\rm final}]$ satisfies 
\begin{equation}
m\Big(\frac{\partial^2 \ufe(t;\mu)}{\partial t^2},v;\mu\Big)+c\Big(\frac{\partial \ufe(t;\mu)}{\partial t},v;\mu\Big)+a\Big(\ufe(t;\mu),v;\mu\Big)=
f(v,t;\mu) \ , \forall v\in X_h \ , \forall t\in[0,T_{\rm final}]\ , \label{eq:TimeVariationalFEGlobal}
\end{equation}
\begin{equation}
\ufe(t=0;\mu)=0 \ ; \frac{\partial \ufe}{\partial t}\Big(t=0;\mu\Big)=0 \ . 
\end{equation}

\subsubsection{Frequency-Domain equation}
\label{sssec:freq}

Similarly to the time-domain equation, the FE approximation $\ufef(\mufreq)$ to $\ufreq(\mufreq)$ can be obtained by projecting equation (\ref{eq:VarFormGlobHom}) on $X_h$: $\ufef(\mufreq)\in \hat X_{h}$ satisfies
\begin{equation}
\afreq(\ufef(\mufreq),v;\mufreq)=\ffreq(v;\mufreq) \  , \forall v\in \hat X_h \ .\label{eq:FreqVariationalFEGlobal}
\end{equation}

\noindent We now present the discrete equations. Let $\Ahattens\in\RR^{\calN_h\times\calN_h}$ and $\underline{\ffreq_h}\in\RR^{\calN_h}$ be the frequency-domain FE matrix and vector defined as 
\begin{equation}
\Big(\Ahattens(\mufreq)\Big)_{qq'}\equiv \afreq(\varphi_{q'},\varphi_q;\mufreq) \ , 1\leq q,q'\leq\calN_h\ , \label{eq:AMatFreq}
\end{equation}
\begin{equation}
\Big(\underline{\ffreq_h}(\mufreq)\Big)_{q}\equiv \ffreq(\varphi_{q};\mufreq) \ , 1\leq q\leq\calN_h\ \label{eq:FvectFreq}
\end{equation}
\noindent The FE basis vector for $\ufef(\mufreq)\in\RR^{\calN_h}$ is then given by
\begin{equation}
\Ahattens(\mufreq)\;\underline{\ufef}(\mufreq)=\underline{\ffreq_h}(\mufreq) \ .
\end{equation}
\noindent Note that $\Ahattens(\mufreq)$ is typically large but very sparse.

\subsection{Finite Element -- Finite Difference Discretization}
\label{subsec:fefddisc}

In order to solve equation (\ref{eq:TimeVariationalFEGlobal}), a finite-difference discretization scheme for time marching with $N_t$ time steps is considered. Let $\Delta t = T_{\rm final}/N_t$ ; $t^j\equiv j\; \Delta t, 0\leq j\leq N_t$. Note that the proposed method can be applied to any finite-difference scheme. A particular scheme is selected in this work for sake of clarity: the unconditionally stable Newmark-$\beta$ scheme with $\beta_t=\frac{1}{4}$ and $\gamma_t=\frac{1}{2}$ (such that the average constant acceleration scheme, or mid-point rule, is obtained) \cite{ChibaKako,Nickelletal}.

Let $\ufefd^j(\mu)$, $0\leq j\leq N_t$ denote the finite element -- finite difference solution at time step $t^j$, and $\dufefd^j(\mu)$ and $\ddufefd^j(\mu)$ the corresponding first and second derivatives in time respectively. Since $\ufefd^0(\mu)=0$ and $\dufefd^0(\mu)=0$, then $\ddufefd^0$ is determined as the solution to
\begin{equation}
m\Big(\ddufefd^0(\mu),v;\mu\Big)=f(v,t=0;\mu) \ , \forall v\in X_h \ .\label{eq:FEFDVariationalddu0}
\end{equation}

\noindent The fields $\ddufefd^j(\mu)$, $\dufefd^j(\mu)$ and $\ufefd^j(\mu)$, $1\leq j\leq N_t$, are then determined as the solutions to the following equations, respectively:
\begin{multline}
m\Big(\ddufefd^j(\mu),v;\mu\Big)+\Delta t\gamma_t \ c\Big(\ddufefd^j(\mu),v;\mu\Big) + \Delta t^2\beta_t \ a\Big(\ddufefd^j(\mu),v;\mu\Big)= f(v,t^j;\mu) \\
-c\Big(\dufefd^{j-1}(\mu)+\Delta t(1-\gamma_t) \ \ddufefd^{j-1}(\mu),v;\mu\Big)-a\Big(\ufefd^{j-1}(\mu)+\Delta t \ \dufefd^{j-1}(\mu)+\Delta t^2(1-\beta_t) \ \ddufefd^{j-1}(\mu),v;\mu\Big) \ , \forall v\in X_h \ ,
\end{multline}
\begin{equation}
\dufefd^j(\mu)=\dufefd^{j-1}(\mu)+\Delta t \ \Big[ (1-\gamma_t) \ \ddufefd^{j-1}(\mu)+\gamma_t \ \ddufefd^j(\mu)\Big]
\end{equation}
\begin{equation}
\ufefd^j(\mu)=\ufefd^{j-1}(\mu)+\Delta t \ \dufefd^{j-1}(\mu)+\Delta t^2 \ \Big[\Big(\frac{1}{2}-\beta_t\Big) \ \ddufefd^{j-1}(\mu)+\beta_t \ \ddufefd^j(\mu) \Big] \ . \label{eq:FEFDVariationaluj}
\end{equation}
 
\noindent We note that the scheme is implicit. 

\noindent We next present the discrete equations in matrix form. Let $\Mtens\in\RR^{\calN_h\times\calN_h}$, $\Ctens\in\RR^{\calN_h\times\calN_h}$ and $\Atens\in\RR^{\calN_h\times\calN_h}$ be the mass, damping and stiffness FE matrices, respectively: 
\begin{equation}
\Big(\Mtens(\mu)\Big)_{qq'}\equiv m(\varphi_{q'},\varphi_q;\mu) \ , 1\leq q,q'\leq\calN_h\ , \label{eq:MMat}
\end{equation}
\begin{equation}
\Big(\Ctens(\mu)\Big)_{qq'}\equiv c(\varphi_{q'},\varphi_q;\mu) \ , 1\leq q,q'\leq\calN_h\ , \label{eq:CMat}
\end{equation}
\begin{equation}
\Big(\Atens(\mu)\Big)_{qq'}\equiv a(\varphi_{q'},\varphi_q;\mu) \ , 1\leq q,q'\leq\calN_h\ . \label{eq:AMat}
\end{equation}

\noindent Furthermore, let $\underline{f_{h}^j}\in\RR^{\calN_h}$, $0\leq j\leq N_t$, be the FE vectors corresponding to the linear form of the time-domain variational formulation at time instance $t^j$, and let $\Ttens\in\RR^{\calN_h\times\calN_h}$ be the time marching matrix:
\begin{equation}
\Big(\underline{f_{h}^j}(\mu)\Big)_{q}\equiv f(\varphi_{q},t_j;\mufreq) \ , 1\leq q\leq\calN_h\ ;\label{eq:Fvect}
\end{equation}
\begin{equation}
\Ttens(\mu)=\Mtens(\mu)+\Delta t\;\gamma_t\;\Ctens(\mu)+\Delta t^2\;\beta_t\;\Atens(\mu)\ .
\end{equation}

\noindent Finally, if $\underline{\ufefd^j}(\mu)\in\RR^{\calN_{h}}$, $\underline{\dufefd^j}(\mu)\in\RR^{\calN_{h}}$ and $\underline{\ddufefd^j}(\mu)\in\RR^{\calN_{h}}$, $0\leq j\leq N_t$, denote the FE basis vectors for ${\ufefd^j(\mu)}$, ${\dufefd^j(\mu)}$ and ${\ddufefd^j(\mu)}$ respectively, we initialize $\underline{\ufefd^0}(\mu)=\underline{0}$, $\underline{\dufefd^0}(\mu)=\underline{0}$, and $\underline{\dufefd^0}(\mu)$ solution of
\begin{equation}
\Mtens(\mu)\;\underline{\ddufefd^0}(\mu) =\underline{f_{h}^0}(\mu) \ ;
\end{equation}
\noindent we then solve for $\underline{\ufefd^j}(\mu)$, $\underline{\dufefd^j}(\mu)$, and $\underline{\ddufefd^j}(\mu)$, for $1\leq j\leq N_t$, from
\begin{equation}
\Ttens(\mu)\;\underline{\ddufefd^j}(\mu) = \Big[ \underline{f_{h}^j}(\mu) - \Ctens(\mu)\; \Big( \underline{\dufefd^{j-1}}(\mu)+\Delta t(1-\gamma_t) \ \underline{\ddufefd^{j-1}}(\mu) \Big) - \Atens(\mu)\;\Big( \underline{\ufefd^{j-1}}(\mu)+\Delta t \ \underline{\dufefd^{j-1}}(\mu) + \Delta t^2(1-\beta_t) \ \underline{\ddufefd^{j-1}}(\mu) \Big) \Big]\ , 
\end{equation}
\begin{equation}
\underline{\dufefd^{j}}(\mu)=\underline{\dufefd^{j-1}}(\mu)+\Delta t \ \Big[(1-\gamma_t) \ \underline{\ddufefd^{j-1}}(\mu)+\gamma_t \ \underline{\ddufefd^{j}}(\mu) \Big]\ , 
\end{equation}
\begin{equation}
\underline{\ufefd^{j}}(\mu)=\underline{\ufefd^{j-1}}(\mu)+\Delta t \ \underline{\dufefd^{j-1}}(\mu)+\Delta t^2 \ \Big[ \Big(\frac{1}{2}-\beta_t\Big) \underline{\ddufefd^{j-1}}(\mu) + \beta_t \ \underline{\ddufefd^{j}}(\mu) \Big]\ .
\end{equation}

\noindent This completes the FE ``truth" discretization.



\section{PR-RBC Approach: Overview}
\label{sec:PRRBC}

In this section, we provide a very brief summary of the Port-Reduced Reduced-Basis Component (PR-RBC) Method. This method will be used for Level $1$ reduction of the two-level reduced basis method proposed in this work and described in Section \ref{sec:2RB}. A complete description can be found in \cite{BhouriThesis}. Since the Level $1$ reduction is carried out on the frequency-domain equation, the parameter spaces considered in this section correspond to the parameter spaces introduced for the time-domain PDE augmented with the angular frequency. Referring to the notations introduced in the previous section, the global parameter considered for the Level $1$ reduction is $\mufreq=(\mu,\omega)\in\Pfreq$, where $\Pfreq\in\RR^{{n}_P+1}$ refers to the augmented compact parameter set.

\subsection{Components -- Ports -- System Assembly}

The PR-RBC method is a domain decomposition technique in which the global system is decomposed into smaller components, which will be referred to as instantiated components. This decomposition creates an ensemble of parameterized instantiated components which can be mapped to an ensemble of parameterized archetype components; multiple instantiated components of the global system can correspond to the same archetype component. Moreover, the domain decomposition creates an ensemble of ports, defined as the intersection of the closures of each two adjacent instantiated components, with the latter forming a parameterized bi-component system. These ports can also be mapped to an ensemble of reference ports associated with archetype bi-component systems.

\begin{figure}[htb] 
	\begin{center}
	\includegraphics[scale=0.23]{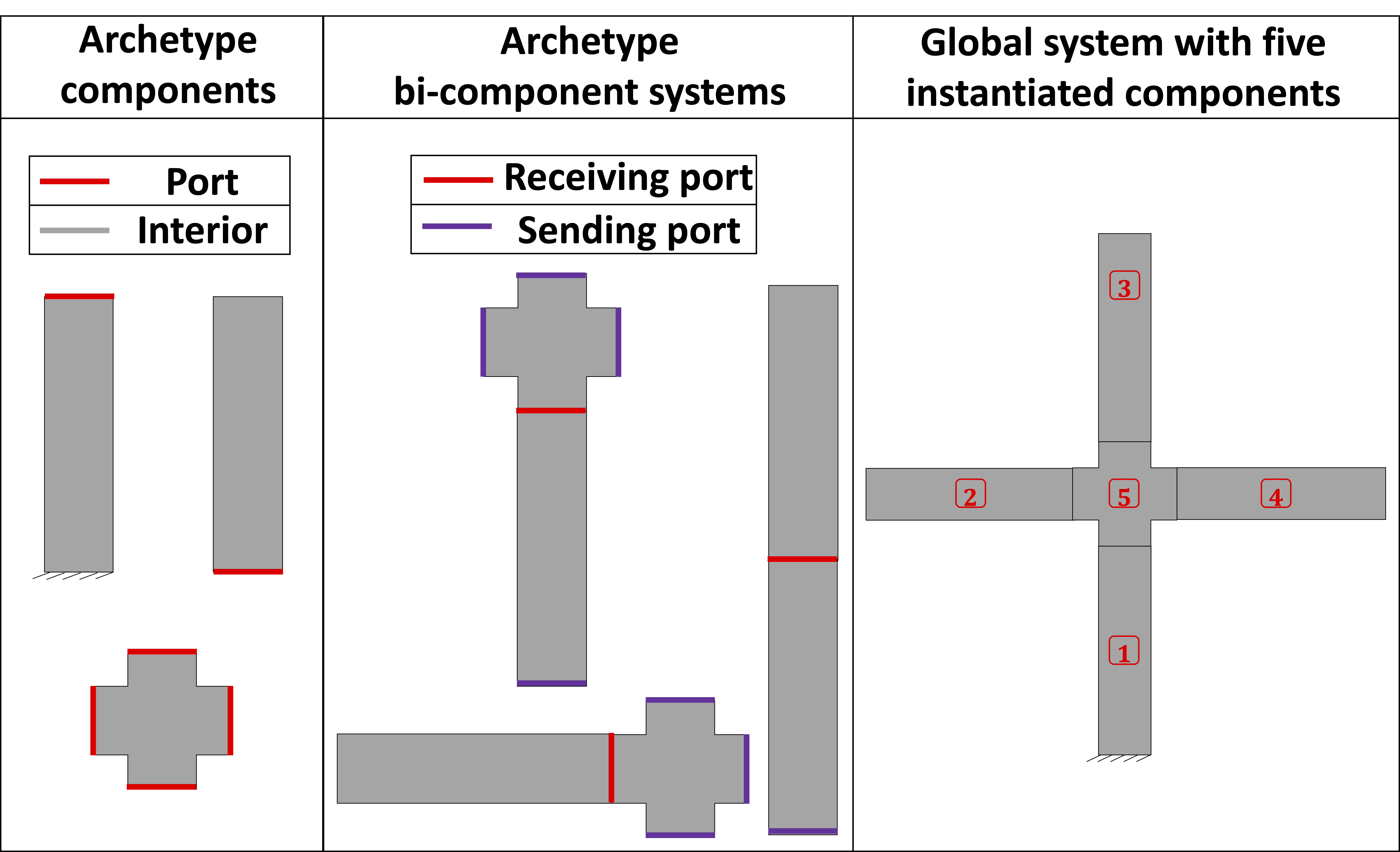}
	\caption{PR-RBC Component--Port--System Assembly: an example with three archetype components, three bi-component systems of compatible archetype components, and a global system with five instantiated components}\label{fig:PR_RBC_overview}
	\end{center}	
\end{figure}

This decomposition technique can be applied to a variety of global systems related through a common physical discipline and hence PDE operator. Therefore, the PR-RBC method can be presented by starting from a library of archetype components and reference ports. Moreover, for simplicity, the ports are presumed to be mutually disjoint, such that the reference port is associated to two local ports. Figure \ref{fig:PR_RBC_overview} shows an example of archetype components (left figure), reference ports (middle figure), and a global system assembly (right figure). In our context, the port and component parameterizations correspond to the frequency-domain equation, and the corresponding sesquilinear and anti-linear forms are assumed to have an affine dependence on the parameters. By consequence, any global system built from the library of the archetype components will be governed by a frequency-domain PDE whose sesquilinear and anti-linear forms have an affine dependence on the system parameter. The latter is intrinsically related to the parameters considered for the different instantiated components forming the global system.

\subsection{PR-RBC Offline Stage}
\label{subsec:PR_RBC_offline}

The PR-RBC offline stage corresponds to the construction of reduced bases to approximate the solution within the parameterized archetype components and over reference ports. Hence, it is informed by the library of archetype components and reference ports, and is independent of any subsequent (feasible) system assembly. Since the system is assembled in the online stage, we need to consider a reference port for each possible bi-component system of compatible archetype components.

For every reference port, a reduced port space needs to be built by port training in order to approximate the solution on the reference port joining each compatible pair of archetype components (red boundary in the middle figure of Figure \ref{fig:PR_RBC_overview}). This construction is carried out by solving a transfer eigenvalue problem such that the reduced port space is optimal in the sense of Kolmogorov $n$-width as shown in \cite{SmetanaPatera}. In addition to the port modes obtained via the transfer eigenvalue problem, additional port modes are considered in order to also account for the inhomogeneity. 

The PR-RBC offline stage also includes the construction of reduced bubble spaces for port mode liftings, and a reduced bubble space for each archetype component with non-zero linear form \cite{EftangPatera2}. These reduced bases approximate the solution inside the archetype components domains (zero on the ports). Within our context, the archetype component bubble spaces reduce to the bubble spaces for inhomogeneity associated with non-zero source terms. In this work, all reduced bubble spaces (for port mode liftings and for inhomogeneity) are constructed by Proper Orthogonal Decomposition (POD). In a more general setting, non-homogeneous Dirichlet boundary conditions for archetype components can be treated by constructing corresponding reduced spaces to approximate the lifting functions \cite{BhouriThesis}. 

Finally, in the context of an offline-online decomposition, all parameter-independent sesquilinear and anti-linear forms needed for the PR-RBC online stage (detailed in Section \ref{subsec:PR_RBC_online}) are computed and stored once in an offline stage. 


\subsection{PR-RBC Online Stage}
\label{subsec:PR_RBC_online} 

For the online stage, we consider a global system characterized by a global parameter and defined as an assembly of instantiated archetype components. 

\begin{figure}[htb] 
	\begin{center}
	\includegraphics[scale=0.28]{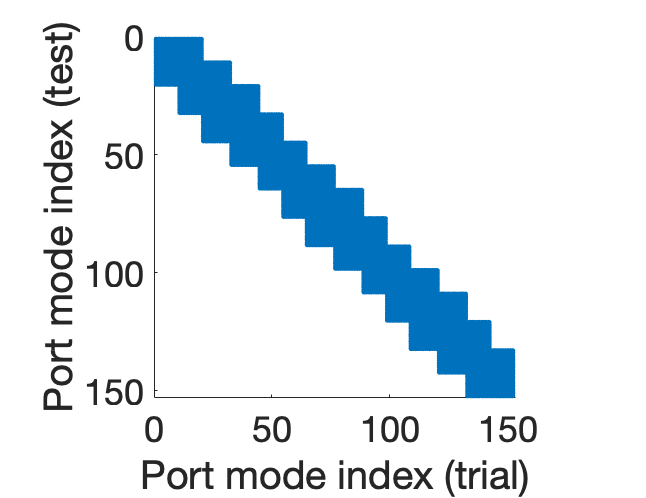}
	\caption{Sparsity pattern of the system Schur complement: $20.6\%$ of the entries are non-zero}
	\label{fig:Schur_spy}
	\end{center}	
\end{figure}

For a given global parameter $\mufreq$, the instantiated component and bi-component parameters are well defined. Using the pre-computed and stored sesquilinear and anti-linear form evaluations, the reduced bubble functions for port mode liftings and for inhomogeneity can be computed at very small computational cost compared to a full FE evaluation. Invoking again the pre-computed and stored sesquilinear and anti-linear forms evaluations, the system Schur complement can be formed and solved. We opt for a Petrov-Galerkin projection to construct the reduced system Schur complement, such that our test space is spanned by the lifted port modes, and not by the ``harmonic" functions which form the (statically condensed) trial space. The Petrov-Galerkin approach is more efficient not only computationally but also in terms of memory storage compared to a Galerkin projection; there is it typically no degradation of stability \cite{Brunken}. Figure \ref{fig:Schur_spy} shows the sparsity pattern of the system PR-RBC Schur complement obtained for the global system considered in Section \ref{sec:num_ex}, resulting from the Petrov-Galerkin projection of the frequency-domain sesquilinear form (\ref{eq:VarFormGlobHom}) using the test space and trial space detailed above. Only port modes which share a common component will result in overlapping reduced spaces, and we thus obtain the ``staircase" sparsity pattern observed in figure \ref{fig:Schur_spy}.

\begin{figure}[htb] 
	\begin{center}
	\includegraphics[scale=0.28]{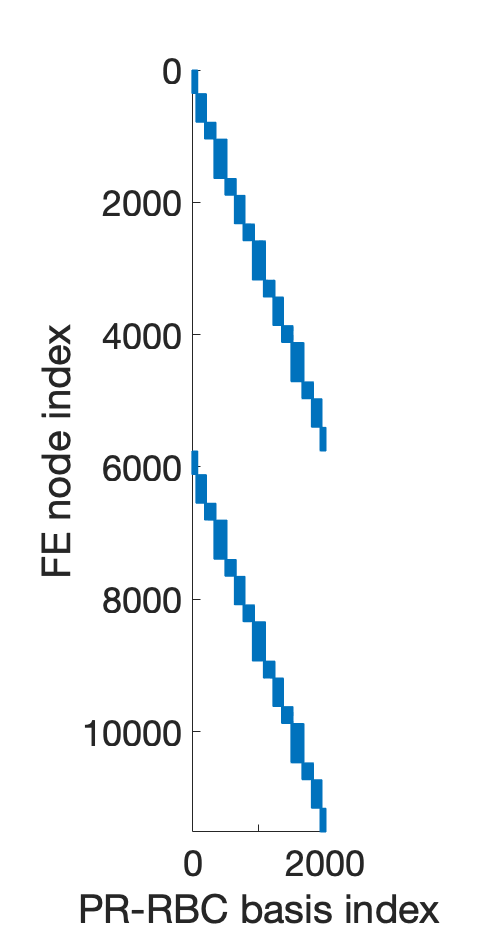}
	\caption{$\underline{\ZZ}$ matrix: $6.9\%$ of the entries are non-zero}
	\label{fig:ZZ_spy}
	\end{center}	
\end{figure}

These solutions define the PR-RBC approximations. A first crucial point related to efficiency is the relatively low dimension of the PR-RBC space. A second crucial point related to efficiency is the sparsity of the PR-RBC basis: the support of a given basis function does not exceed two instantiated components for lifted port modes, and is further restricted to just one instantiated component for reduced bubble spaces. Let $X_{h,\bD}$, of dimension $N_{h,\bD}$, be the PR-RBC space for the global system considered, where {\bf D} signifies the many port and bubble discretization parameters which inform the PR-RBC approximation. Figure \ref{fig:ZZ_spy} shows the sparsity pattern of the $\underline{\ZZ}\in\RR^{\calN_h\times N_{h,\bD}}$ matrix containing the FE representation of the PR-RBC basis obtained for the global system detailed in Section \ref{sec:num_ex}. The columns of $\underline{\ZZ}$ correspond to the coefficients of the reduced port modes and reduced bubble modes used for the global domain system considered in Section \ref{sec:num_ex} as represented by the FE nodal basis $\{\varphi_j\}_{j=1,\ldots,\calN_h}$: a member $\upr(\mufreq)\in X_{h,\bD}$ can be expressed as
\begin{equation}
\upr(\mufreq)=\sum\limits_{j=1}^{\calN_h} \big(\underline{\ZZ}\; \underline{\upr}(\mufreq)\big)_j\varphi_j \ ,
\end{equation}
\noindent where $\underline{\upr}(\mufreq)\in\RR^{N_{h,\bD}}$ refers to the PR-RBC basis vector corresponding to $\upr(\mufreq)$. 

We emphasize the important role of components. In general, the components distribute the parameter domain: we reduce a large problem with many global parameters to many small problems each with just a few (local) parameters. Components also permit consideration of very large systems: even in the PR-RBC offline stage, we are required to solve FE problems over at most pairs of components | never the full system. Also, components provide geometry and topology parametric variation. And finally, components permit us to more easily justify the PR-RBC Offline investment: we may amortize the offline effort not only over many queries for any particular global system, but over all possible global systems in our family. Note in this sense we can formally define our family of feasible global systems as the set of all systems which may be constructed from the associated library of archetype components. This family of global systems can then be well approximated either in the PR-RBC online stage, or in the proposed two-level reduced basis method proposed in this work.



\section{Two-Level Reduced Basis Method}
\label{sec:2RB}

To compute an approximation to the time-domain elastodynamics equation for a given global system and a parameter value $\mu$, we propose a two-level model order reduction approach. The first-level reduction consists of evaluating the PR-RBC solutions to the frequency-domain PDE for the global system at well-selected frequencies, and hence involves the PR-RBC online stage. The second-level reduction consists of building a ``final" reduced basis by a Strong Greedy approach in which the inexpensive PR-RBC solutions to the frequency-domain equation serve as surrogate ``truth" solutions. The time-domain elastodynamics PDE is then projected and solved within the ``final" reduced basis in standard fashion. Note that the PR-RBC offline stage is conducted prior to execution of the two-level procedure. Then, in the online stage, for the given parameter value $\mu$, both levels of reduction are invoked, and hence both levels must be computationally fast. In contrast, the PR-RBC offline stage is run only once independently of the number of parameter values considered in the evaluation of the time-domain solution, hence, we do not give as much importance to the computational cost of the PR-RBC offline stage.

\subsection{Level 1}

For a given global system test parameter $\mu$, we consider a sufficiently rich angular frequency $\omega$ set $\Xi_\omega$, of size $n_\omega$, and we consider the online-train (o-t) dataset $\Xi_{\rm o-t}=\{\mufreq\equiv(\mu,\omega) ; \omega\in\Xi_\omega\}$. Hence $\Xi_{\rm o-t}$ can be expressed as: $\Xi_{\rm o-t}=\{\mufreq_j, 1\leq j\leq n_\omega  \}$, where $\mufreq_j=(\mu,\omega_j)$ for $\omega_j\in\Xi_\omega$. Let $\sigma_t^{\rm ref}$ refer to a characteristic time of $f_t(\cdot;\mu_t)$ (the latter being defined in equation (\ref{eq:fxft})); then $\Xi_\omega$ is chosen as
\begin{equation}\label{eq:xi_omega}
\Xi_\omega=\{0,d\omega,\ldots,\omega_{\rm max}\} \ , \ d\omega=\frac{1}{\underline{c_\omega}\; {\sigma_t^{\rm ref}}} \ , \ \omega_{\rm max}=\frac{\overline{c_\omega}}{{\sigma_t^{\rm ref}}}\ , 
\end{equation}

\noindent so that $n_\omega=\overline{c_\omega} \; \underline{c_\omega}+1$.

Note that for the PR-RBC offline stage described in Section \ref{subsec:PR_RBC_offline}, the training is performed over a frequency set that at least contains $\Xi_\omega$. Nonetheless, the PR-RBC approach is still general in the sense that we can consider any assembly of components forming a feasible global system as mentioned in section \ref{subsec:PR_RBC_online}. For the numerical examples considered in this work, the PR-RBC offline stage is performed over a frequency set that exactly matches $\Xi_\omega$. Then, the PR-RBC online stage (Section \ref{subsec:PR_RBC_offline}) is performed to obtain approximations to the frequency-domain equation for $\mufreq_j\in\Xi_{\rm o-t}$, which gives $n_\omega$ approximations, $\upr(\mufreq_j)$ for $\mufreq_j\in\Xi_{\rm o-t}$, to $n_\omega$ FE solutions $\ufef(\mufreq_j)$ for $\mufreq_j\in\Xi_{\rm o-t}$. This step corresponds to the first-level of reduction.

\subsection{Level 2}

The second-level reduction consists of constructing a reduced basis from the PR-RBC approximations computed in Level 1 by performing a Strong Greedy procedure to identify a reduced space $X_{\rm RB}$ of size $N$, where $\XRB\in\RR^{\calN_{h}\times N}$ is the FE representation of $X_{\rm RB}$. The column $j$ of $\XRB$ corresponds to the coefficients of the RB basis function $j$ from $X_{\rm RB}$ as represented by the FE nodal basis $\{\varphi_i\}_{i=1,\ldots,\calN_h}$. Note that in this Level 2 reduction, only the angular frequency $\omega$ varies. Therefore, the final reduced space $X_{\rm RB}$ is expected to be of sufficiently small size and also able to provide a sufficiently accurate reduced basis approximation to the full FE solution. Finally, the time-domain PDE is projected on $X_{\rm RB}$. The exact time signature of the linear form is only used in the time marching performed using the reduced space $X_{\rm RB}$.

The construction of $X_{\rm RB}$ is carried out using a Strong Greedy approach on $\{\upr(\mufreq_j) \ ,\mufreq_j\in\Xi_{\rm o-t}\}$ evaluated at the first-level reduction. These inexpensive PR-RBC approximations are considered as ``truth" solutions. Let $X_{\rm RB\,i}$ denote the reduced space of size $i$, constructed prior to the $i$-th iteration of the Strong Greedy algorithm; let $\upri(\mufreq_j)\in X_{\rm RB\,i}$ denote the RB approximation to $\upr(\mufreq_j)$ obtained using the reduced space $X_{\rm RB\,i}$. Since the PR-RBC approximations are considered as high-fidelity solutions within the greedy algorithm, the next snapshot we add to the reduced space $X_{\rm RB\,i}$ among the PR-RBC solutions $\{\upr(\mufreq_j) \ ,\mufreq_j\in\Xi_{\rm o-t}\}$ is selected based on the norm of the error: $\big|\big|\upr(\mufreq)-\upri(\mufreq)\big|\big|_{H^1(\Omega)}$, for $\mufreq\in\Xi_{\rm o-t}$. The Strong Greedy approach is summarized in Algorithm \ref{alg:SG}. An efficient computation of the errors $\big|\big|\upr(\mufreq)-\upri(\mufreq)\big|\big|_{H^1(\Omega)}$, for $\mufreq\in\Xi_{\rm o-t}$, such that only the updated quantities depending on the new snapshot are evaluated at every iteration of the Strong Greedy, is detailed in \cite{BhouriThesis} . 

\begin{algorithm}[H]
\caption{Strong Greedy algorithm}\label{alg:SG}
\begin{algorithmic}[1]
\Function{Strong Greedy}{$\{\upr(\mufreq_j) \ ,\mufreq_j\in\Xi_{\rm o-t}\},n_\omega,M,\epsilon
$}

\State $i\leftarrow 1$

\State $e_1\leftarrow \infty$

\State $j$ is a random integer sampled from $\{1,\ldots,n_\omega\}$

\State $X_{\rm RB\,1}\leftarrow {\rm span}(\upr(\mufreq_j))$

\State $\epsilon_a\leftarrow \max\limits_{\mufreq\in\Xi_{\rm o-t}}\big|\big|\upr(\mufreq)-\upro(\mufreq)\big|\big|_{H^1(\Omega)}$

\While {$i\leq\min(n_\omega,M)$ and $e_i>\epsilon$}

\State Choose $j\in\{1,\ldots,n_\omega\}$ such that, ${\big|\big|\upr(\mufreq_j)-\upri(\mufreq_j)\big|\big|_{H^1(\Omega)}}={\max\limits_{\mufreq\in\Xi_{\rm o-t}}\big|\big|\upr(\mufreq)-\upri(\mufreq)\big|\big|_{H^1(\Omega)}}$

\State $X_{\rm RB\,i+1}\leftarrow X_{\rm RB\,i}\oplus {\rm span}( \upr(\mufreq_j))$

\State $e_{i+1}\leftarrow \frac{\big|\big|\upr(\mufreq_j)-\upri(\mufreq_j)\big|\big|_{H^1(\Omega)}}{\epsilon_a}$

\State $i\leftarrow i+1$

\EndWhile
\State \textbf{end while}


\State \textbf{return} $X_{\rm RB\,i}$
\\ \textbf{end function}
\EndFunction
\end{algorithmic}
\end{algorithm}

We can compute the errors $\big|\big|\upr(\mufreq)-\upri(\mufreq)\big|\big|_{H^1(\Omega)}$ for $\mufreq\in\Xi_{\rm o-t}$ within the strong greedy iterations using pre-computed and stored parameter-independent inner products formed from the PR-RBC basis. The resulting operation count would be independent of $\calN_h$, but dependent of the number of PR-RBC basis functions that have overlapping supports. We found it was considerably less expensive, and simpler to use the full FE representation since we only need sparse matrix-vector multiplication in order to determine $\big|\big|\upr(\mufreq)-\upri(\mufreq)\big|\big|_{H^1(\Omega)}$, for $\mufreq\in\Xi_{\rm o-t}$. Hence, the errors computation is carried out using the FE representations given by $\underline{\upr^{\rm FE}}(\mufreq)=\underline{\ZZ}\; \underline{\upr}(\mufreq)$ , for $\mufreq\in\Xi_{\rm o-t}$; here $\underline{\upr}(\mufreq)\in\RR^{N_{h,\bD}}$ refers to the basis vector for $\upr(\mufreq)$, and $\underline{\ZZ}\in\RR^{\calN_h\times N_{h,\bD}}$ is the sparse matrix containing the FE representation of the PR-RBC basis as introduced in section \ref{subsec:PR_RBC_online} (see figure \ref{fig:ZZ_spy}). 
 
Similarly to the FE approximation $\ufe(\cdot,\mu)$, the two-level reduced basis approximation, denoted $\urb(t\in[0,T_{\rm final}];\mu)$, is obtained by projecting equation (\ref{eq:VarFormGlobElast}) on $X_{\rm RB}$: $\urb(t;\mu)\in X_{\rm RB}\ , \forall t\in[0,T_{\rm final}]$ ,  such that
\begin{equation}
m\Big(\frac{\partial^2 \urb(t;\mu)}{\partial t^2},v;\mu\Big)+c\Big(\frac{\partial \urb(t;\mu)}{\partial t},v;\mu\Big)+a\Big(\urb(t;\mu),v;\mu\Big)=
f(v,t;\mu) \ , \forall v\in X_{\rm RB} \ , \forall t\in[0,T_{\rm final}]\ , \label{eq:TimeVariationalRBGlobal}
\end{equation}
\begin{equation}
\urb(t=0;\mu)=0 \ ; \frac{\partial \urb}{\partial t}\Big(t=0;\mu\Big)=0 \ .
\end{equation}

\noindent We now incorporate a finite-difference scheme with the same time-discretization notations introduced in Section \ref{subsec:fefddisc}. Let $\urbfd^j(\mu)$ for $0\leq j\leq N_t$ denote the two-level reduced basis-finite difference solution at time step $t^j$, and $\durbfd^j$ and $\ddurbfd^j$ the corresponding first and second derivatives in time respectively. We initialize $\urbfd^0(\mu)=0$, $\durbfd^0(\mu)=0$, and $\ddurbfd^0$ solution of
\begin{equation}
m\Big(\ddurbfd^0(\mu),v;\mu\Big)= f(v,t=0;\mu) \ , \forall v\in X_{\rm RB} \ ;\label{eq:RBFDVariationalddu0}
\end{equation}
\noindent we then solve for $\ddurbfd^j(\mu)$, $\durbfd^j(\mu)$, and $\urbfd^j(\mu)$, for $1\leq j\leq N_t$, from
\begin{multline}
m\Big(\ddurbfd^j(\mu),v;\mu\Big)+\Delta t\gamma_t \ c\Big(\ddurbfd^j(\mu),v;\mu\Big) + \Delta t^2\beta_t \ a\Big(\ddurbfd^j(\mu),v;\mu\Big)=f(v,t^j;\mu) \\
-c\Big(\durbfd^{j-1}(\mu)+\Delta t(1-\gamma_t) \ \ddurbfd^{j-1}(\mu),v;\mu\Big)-a\Big(\urbfd^{j-1}(\mu)+\Delta t \ \durbfd^{j-1}(\mu)+\Delta t^2(1-\beta_t) \ \ddurbfd^{j-1}(\mu),v;\mu\Big) \ , \forall v\in X_{\rm RB} \ ,
\end{multline}
\begin{equation}
\durbfd^j(\mu)=\durbfd^{j-1}(\mu)+\Delta t \ \Big[ (1-\gamma_t) \ \ddurbfd^{j-1}(\mu)+\gamma_t \ \ddurbfd^j(\mu)\Big] \ ,
\end{equation}
\begin{equation}
\urbfd^j(\mu)=\urbfd^{j-1}(\mu)+\Delta t \ \durbfd^{j-1}(\mu)+\Delta t^2 \ \Big[\Big(\frac{1}{2}-\beta_t\Big) \ \ddurbfd^{j-1}(\mu)+\beta_t \ \ddurbfd^j(\mu) \Big] \ , \label{eq:RBFDVariationaluj}
\end{equation}
\noindent respectively.
 
We now provide the matrix equations. Let $\MRB\in\RR^{N\times N}$, $\CRB\in\RR^{N\times N}$ and $\ARB\in\RR^{N\times N}$ be the mass, damping and stiffness two-level reduced basis matrices, respectively, 
\begin{equation}
\MRB(\mu)=  \XRB(\mu)^H\;\Mtens(\mu)\;\XRB(\mu), \label{eq:MRBMat}
\end{equation}
\begin{equation}
\CRB(\mu)=  \XRB(\mu)^H\;\Ctens(\mu)\;\XRB(\mu), \label{eq:CRBMat}
\end{equation}
\begin{equation}
\ARB(\mu)=  \XRB(\mu)^H\;\Atens(\mu)\;\XRB(\mu), \label{eq:ARBMat}
\end{equation}
\noindent where $\cdot^H$ denotes the Hermitian transpose operator. Similarly, let $\underline{f_{\rm RB}^j}\in\RR^{N}$, $0\leq j\leq N_t$, be the two-level reduced basis vectors corresponding to the linear form of the time-domain variational formulation at time instance $t^j$, and $\TRB\in\RR^{N\times N}$ the time marching matrix,
\begin{equation}
\underline{f_{\rm RB}^j}(\mu)= \XRB(\mu)^\dagger\; \underline{f_{h}^j}(\mu) \label{eq:FRBvect}
\end{equation}
\begin{equation}
\TRB(\mu)=\MRB(\mu)+\Delta t\;\gamma_t\;\CRB(\mu)+\Delta t^2\;\beta_t\;\ARB(\mu)\ ,
\end{equation}
\noindent respectively. 

Let $\underline{\urbfd^j}(\mu)\in\RR^{N}$, $\underline{\durbfd^j}(\mu)\in\RR^{N}$ and $\underline{\ddurbfd^j}(\mu)\in\RR^{N}$, $0\leq j\leq N_t$, denote the two-level reduced basis vectors of ${\urbfd^j}(\mu)$, ${\durbfd^j}(\mu)$ and ${\ddurbfd^j}(\mu)$ respectively. We initialize $\underline{\urbfd^0}(\mu)=\underline{0}$, $\underline{\durbfd^0}(\mu)=\underline{0}$, and $\underline{\ddurbfd^0}(\mu)$ solution of
\begin{equation}
\MRB(\mu)\;\underline{\ddurbfd^0}(\mu) = \underline{f_{\rm RB}^0}(\mu) \ ;
\end{equation}
\noindent we then proceed to time-march as
\begin{multline}
\TRB(\mu)\;\underline{\ddurbfd^j}(\mu) =\\
\Big[ \underline{f_{\rm RB}^j}(\mu) - \CRB(\mu)\; \Big( \underline{\durbfd^{j-1}}(\mu)+\Delta t(1-\gamma_t) \ \underline{\ddurbfd^{j-1}}(\mu) \Big) - \ARB(\mu)\;\Big( \underline{\urbfd^{j-1}}(\mu)+\Delta t \ \underline{\durbfd^{j-1}}(\mu) + \Delta t^2(1-\beta_t) \ \underline{\ddurbfd^{j-1}}(\mu) \Big) \Big]\ , 
\end{multline}
\begin{equation}
\underline{\durbfd^{j}}(\mu)=\underline{\durbfd^{j-1}}(\mu)+\Delta t \ \Big[(1-\gamma_t) \ \underline{\ddurbfd^{j-1}}(\mu)+\gamma_t \ \underline{\ddurbfd^{j}}(\mu) \Big]\ , 
\end{equation}
\begin{equation}
\underline{\urbfd^{j}}(\mu)=\underline{\urbfd^{j-1}}(\mu)+\Delta t \ \underline{\durbfd^{j-1}}(\mu)+\Delta t^2 \ \Big[ \Big(\frac{1}{2}-\beta_t\Big) \underline{\ddurbfd^{j-1}}(\mu) + \beta_t \ \underline{\ddurbfd^{j}}(\mu) \Big]\ , 
\end{equation}
\noindent for $1\leq j\leq N_t$.

Since $X_{\rm RB}$ is a complex-valued reduced basis, the actual two-step reduced basis approximation of the finite element - finite difference solution $\ufefd^j(\mu)$ is given by:
\begin{equation}
\urbfd^j(\mu)=\sum\limits_{k=1}^{\calN_{h}}\Big(\Re\Big[\XRB(\mu)\;\underline{\urbfd^{j}}(\mu)\Big]\Big)_k\varphi_{k}\ , 1\leq j\leq N_t\ .
\end{equation}

\noindent If a quantity of interest $\underline{q}\in\RR^{N_q}$ is considered, then the corresponding FE output matrix $\Qtens\in\RR^{N_q\times\calN_{h}}$ is multiplied by $\XRB(\mu)$ once (offline) to obtain
\begin{equation}
\QRB(\mu)=  \Qtens\;\XRB(\mu), \label{eq:QRBMat}
\end{equation}
\noindent and the two-level reduced basis approximation of $\underline{q}$ at time step $t^j$ is simply given by:
\begin{equation}
\underline{q^j_{\rm RB}}= \Re\Big[ \QRB(\mu)\;\underline{\urbfd^{j}}(\mu)\Big] \ ; 1\leq j\leq N_t. \label{eq:QRB}
\end{equation}

Note that the error between the two-level PR-RBC approximation and the FE solution has three sources. One source is related to the discretization $\Xi_\omega$ and truncation of the inverse Laplace transform. The two other sources relate to reduction: one is introduced in Level 2 by the Strong Greedy-based reduction, and the other is introduced in Level 1 by considering the PR-RBC approximations as high-fidelity truth. This decomposition of the different sources of the error can be useful for developing an {\em a priori} error estimate.



\subsection{Operation Count}
\label{subsec:op_count}

Let ${n}_{\rm comp}$ and ${n}_{\rm port}$ be the number of instantiated components and the number of ports forming a global system. The operation count to estimate $\urbfd^j(\mu)$ for $1\leq j\leq N_t$ and $\ntt$ different parameters $\mu$ by the two-level reduced basis method described above is given by:
\begin{multline}
\underbrace{\ntt\; n_\omega\; \Big[ \sum_{i=1}^{{n}_{\rm comp}}\mathcal{O}\big(({N}_i^{\rm inhom})^3\big) + \sum_{k=1}^{{n}_{\rm port}}\sum_{m=1}^{M_k'}\sum_{\ell=1}^{2} \mathcal{O}\big((N_{k,m,\ell})^3\big) + \Big( \sum_{k=1}^{n_{\rm port}} M_k' \Big)^\kappa\Big]}_{(\rom{1})}+\\
\underbrace{\ntt\; n_\omega\;\mathcal{O}\big(n_z + n_h + N\; \calN_h + N^4\big)}_{(\rom{2})}+ \underbrace{\ntt\; N_t\;\mathcal{O}\big(N^2+N_q\; N\big)}_{(\rom{3})} \ . \label{eq:operation_count_online}
\end{multline}

\noindent Here $M_k'$ refers to the size of reduced port space for port $k$, ${N}_i^{\rm inhom}$ is the size of the reduced bubble space for inhomogeneity for component $i$, and $N_{k,m,\ell}$ , for $\ell=1,2$, refers to the size of reduced bubble spaces for port mode liftings for port $k$ and port mode $m$ within the two components forming the port $k$. The exponent $\kappa$ is a solver-dependent scaling exponent for sparse matrix inversion | in particular, inversion of the sparse PR-RBC Schur complement (see figure \ref{fig:Schur_spy}). Finally, $n_z$ and $n_h$ refer to the number of non-zero entries of the $\underline{\ZZ}$ matrix and the FE norm matrix respecively.


The term $(\rom{1})$ of the operation count given in (\ref{eq:operation_count_online}) corresponds to the first-level reduction in which the PR-RBC online stage is run $n_\omega$ times; the term $(\rom{2})$ corresponds to performing the Strong Greedy algorithm and projecting the time-domain equation into $X_{\rm RB}$; the term $(\rom{3})$ corresponds to (i) time marching using the reduced space $X_{\rm RB}$ via a forward-backward substitution (neglecting the computational cost of the initial LU decomposition), and (ii) computation of the output of interest. Note that the proposed two-level reduced basis method does not contain any FE matrix inversion. In contrast, computing the FE approximations $\ufefd(\cdot;\mu)$ for $\ntt$ different parameters $\mu$ requires:
\begin{equation}
\ntt\; N_t\cdot ((\calN_{h})^{\kappa''}+n_{\rm outputs}\cdot \calN_{h}) \ ,
\end{equation}
\noindent where $\kappa''>1$ is a solver-dependent scaling exponent for (FE) sparse matrix forward-backward substitution.



As shown in the numerical example detailed in section \ref{sec:num_ex}, we generally have $N_t\gg n_\omega$ for large geometric domains with localized excitations with relatively short time signatures. Moreover, unlike the proposed two-level reduced basis method, the operation count of the FE approximation contains the term $(\calN_{h})^{\kappa''}$. Hence, the two-level reduced basis method is expected to have considerably smaller computation cost compared to the full FE approximation, as we will show in the numerical experiments in section \ref{sec:num_ex}.

\section{Numerical Example}
\label{sec:num_ex}


\subsection{Archetype Components and Bi-component Systems}

\begin{figure}[htb]
	\begin{center}
	\includegraphics[scale=0.22]{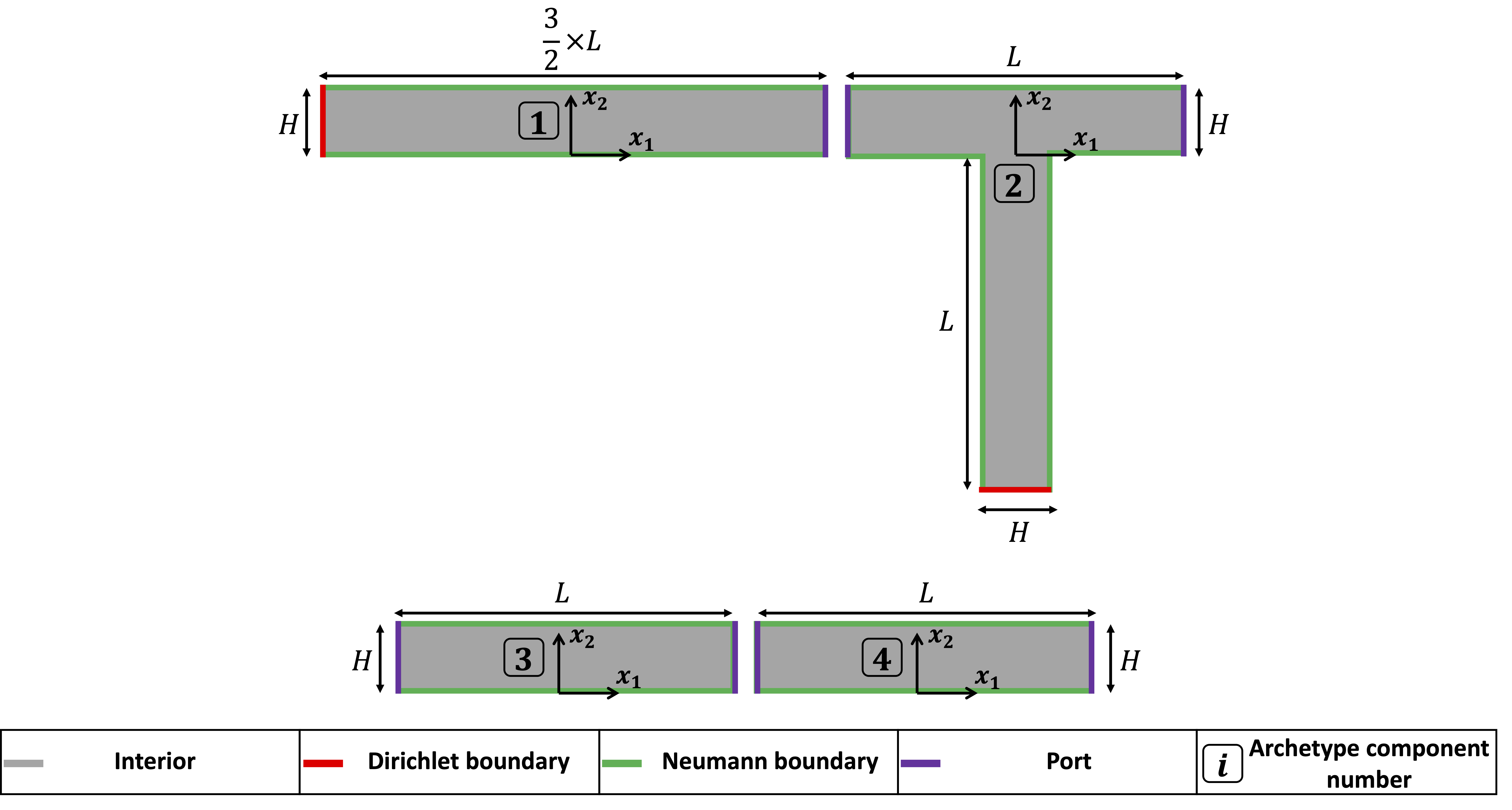}
	\caption{Library of archetype components (considered for construction of elastodynamics bridges)}
	\label{fig:bridge_PR_RBC_ref_comp}
	\end{center}	
\end{figure}

\begin{figure}[htb]
	\begin{center}
	\includegraphics[scale=0.22]{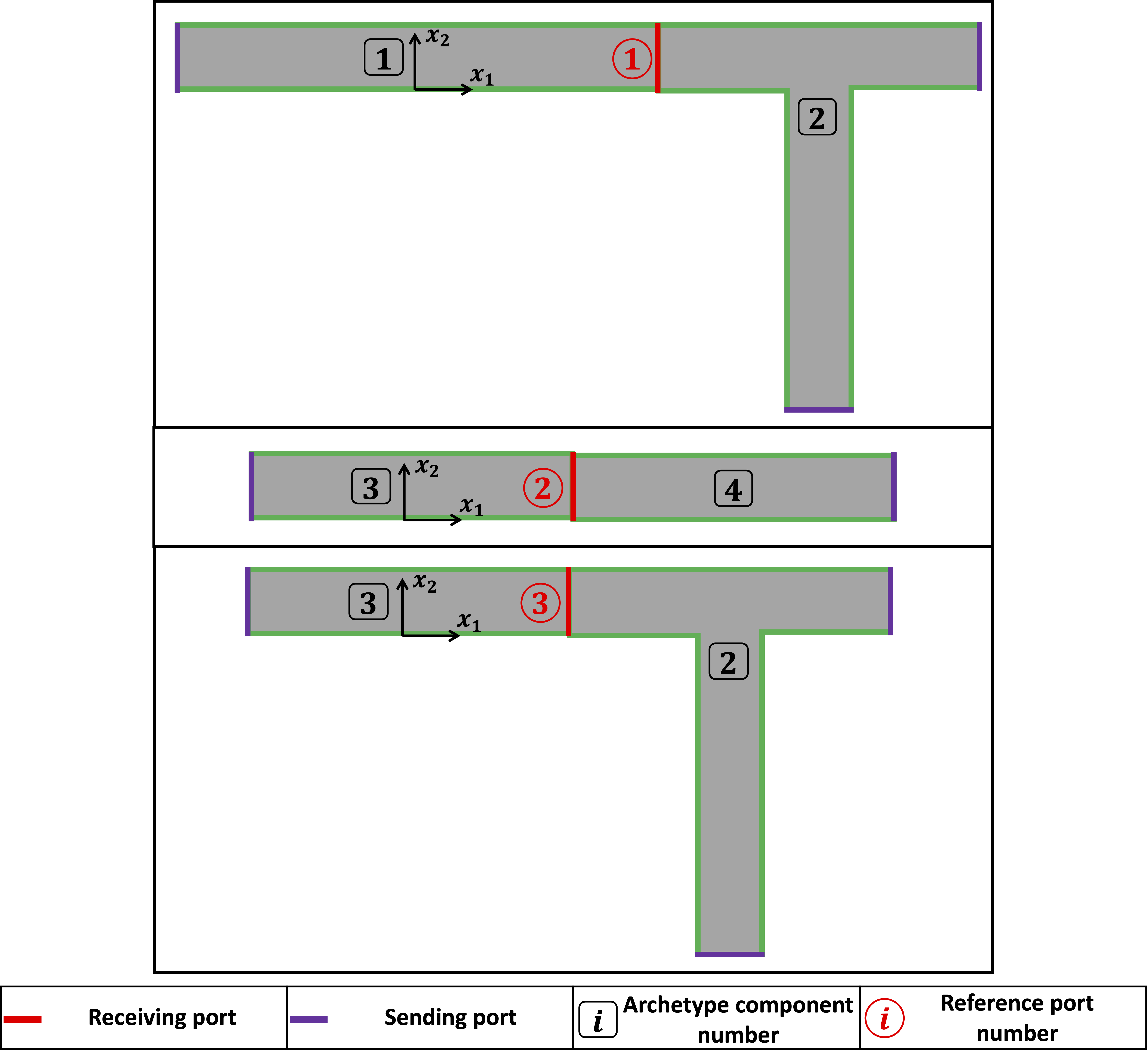}
	\caption{Library of reference ports (archetype bi-components) constructed from the archetype components defined in figure \ref{fig:bridge_PR_RBC_ref_comp}}		
	\label{fig:bridge_PR_RBC_ref_port}
	\end{center}	
\end{figure}

In this section, we apply the two-level PR-RBC approach to the 2D-elastodynamics equation for the library of archetype components shown in figure \ref{fig:bridge_PR_RBC_ref_comp} and the associated reference ports shown in figure \ref{fig:bridge_PR_RBC_ref_port}. There are two different approaches to port training: 1) train ports over all archetype bi-components which have the same reference port geometry, 2) train ports over all distinct archetype bi-components. The former approach is useful if the number of distinct archetype bi-components is relatively very large compared to the allocated memory ressources, while the latter generally provides smaller reduced port spaces for a given accuracy criterion since every reduced port space is constructed to approximate a specific port resulting from a single archetype bi-component. In this work, we opt for the second approach since we have only three distinct archetype bi-components. The archetype component $1$ has a rectangular geometry of dimension $\frac{3}{2} L\times H$ and contains a homogeneous Dirichlet boundary and homogeneous Neumann boundaries. The archetype component number $2$ has a T shape as detailed in figure \ref{fig:bridge_PR_RBC_ref_comp} and also contains a homogeneous Dirichlet boundary and homogeneous Neumann boundaries. The archetype components $3$ and $4$ have a rectangular geometry of dimension $L\times H$ and do not contain any Dirichlet boundary. Components $3$ has only homogeneous Neumann boundary conditions, while component $4$ has a non-homogeneous Neumann term corresponding to traction,
\begin{equation}
\sigma(u)\cdot n=\begin{bmatrix} F\times t\times e^{-\frac{t}{\sigma_t}}\times e^{-\frac{(x_1-x_c)^2}{\sigma_x^2}}\times\mathbbm{1}_{\{x_2=H\}} \\ -c_{\rm friction}\times F\times t\times e^{-\frac{t}{\sigma_t}}\times e^{-\frac{(x_1-x_c)^2}{\sigma_x^2}}\times\mathbbm{1}_{\{x_2=H\}} \end{bmatrix}  \ , \forall x_1\in \Big[-\frac{L}{2},\frac{L}{2}\Big] \ , x_2\in\{0,H\} \ . \label{eq:nh_NBC}
\end{equation}
\noindent Here $\sigma(u)$ denotes the stress tensor, $F$ the load amplitude, $\sigma_t$ the temporal parameter, $x_c$ the load center, $\sigma_x$ the load spatial width, $c_{\rm friction}$ the fraction coefficient, $\mathbbm{1}_C$ the 2D-function equal to $1$ if the condition $C$ is satisfied and $0$ otherwise, and $n$ the outer normal of the geometric domain.


Each of the archetype components $1$, $2$ and $3$ has $3$ parameters which consist of the Young's modulus and the Rayleigh damping coefficients as introduced in equations (\ref{eq:a_bil_form_time_domain}) and (\ref{eq:c_bil_form_time_domain}). Hence, the corresponding frequency-domain variational problem has $4$ parameters | including the angular frequency. Archetype component $4$ has $5$ additional parameters defining the non-homogeneous Neumann boundary condition detailed in equation (\ref{eq:nh_NBC}): the load amplitude $F$, the temporal parameter $\sigma_t$, the load center $x_c$, the load spatial width $\sigma_x$, and the friction coefficient $c_{\rm friction}$. Since the frequency-domain representation of the load is taken constant and equal to $1$ for any angular frequency value, the frequency-domain variational problem for archetype component $4$ has $8$ parameters in total. The reference ports parameters and the boundary conditions considered for the associated bi-component problems follow naturally from the parameters and the boundary conditions defined for the archetype components. The variational problem for archetype component-pairwise training has either $7$ parameters, for reference ports not involving the archetype component $4$, or $11$ parameters for the reference port involving the archetype component $4$. These local parameters associated with the frequency-domain variational problems define the parameters spaces considered in building the reduced bases that are used at Level 1 reduction. 

We choose acrylic as material for which $\overline{\alpha_{\rm Ray}}=5.3785\times10^{-4}s^{-1}$ and $\overline{\beta_{\rm Ray}}=1.0634\times10^{-4}s$ correspond to nominal values of the Rayleigh damping coefficients \cite{TaddeiThesis}. For acrylic, density and Poisson's ratio are well-characterized in the literature (\cite{HarrisSabnis} Chapter 3.6.2 for the Poisson's ratio and the webpage pubchem.ncbi.nlm.nih.gov for the density). We therefore set $\rho=1180\; kg.m^{-3}$, $\bar{E}=2.755\; GPa$, and $\bar\nu=0.35$ as nominal values, and we choose $\sigma_t^{\rm ref}=16\; T_{\rm ref}$ (see equation (\ref{eq:xi_omega})), where $T_{\rm ref}\equiv\frac{H}{c_t}$ and $c_t=\sqrt{\frac{\bar{E}}{2\;\rho\;(1+\bar\nu)}}$ is the celerity of the transverse wave in infinite domain without damping. We take $\underline{c_\omega}=10$ and $\overline{c_\omega}=4$ (see equation (\ref{eq:xi_omega})), such that $n_\omega=41$ and for $\mu_t\equiv\sigma_t\in[0.75\;\sigma_t^{\rm ref},1.25\;\sigma_t^{\rm ref}]$, 
\begin{equation}
|\hat f_t(\omega,\mu_t)|\leq0.06\,\max\limits_{\omega'\geq0}|\hat f_t(\omega',\mu_t)| \, , \omega\geq \omega_{\rm max} \ .
\end{equation}

We take $L=5\; m$, $H=1\; m$ and consider the following domains from which parameters are uniformly sampled for training of archetype components and associated bi-components:
\begin{equation}
{\alpha_{\rm Ray}}\in]0,\overline{\alpha_{\rm Ray}}] \ , 
\end{equation}
\begin{equation}
{\beta_{\rm Ray}}\in]0,\overline{\beta_{\rm Ray}}] \ , 
\end{equation}
\begin{equation}
E\in[0.75\;\bar{E},1.25\;\bar{E}] \ , 
\end{equation}
\begin{equation}
F\in\Big[-20\;\frac{\bar{E}}{T_{\rm ref}},-10\;\frac{\bar{E}}{T_{\rm ref}}\Big] \ , 
\end{equation}
\begin{equation}
\sigma_x\in[0.02\; m,0.04\; m] \ , 
\end{equation}
\begin{equation}
x_c\in[2.46\; m,2.54\; m] \ ,
\end{equation}
\begin{equation}
c_{\rm friction}\in[0.5,0.7] \ .
\end{equation} 

\noindent We assume that we have non-zero damping such that we have the theoretical foundation justifying the existence of the Laplace transform. Nonetheless, the proposed two-level reduced basis method still performs well even for the zero damping case, mainly due to the stability introduced by the finite difference scheme in time. Hence, the proposed two-level reduced basis method can be applied to non-damped or damped systems; we consider small damping.

The sizes of the different reduced bases formed at Level 1 reduction are chosen based on the decrease of the eigenvalues of the transfer eigenvalue problem and the decrease of the POD modes for the reduced bubble space for inhomogeneity, the reduced port space, and the reduced space for port mode lifting. Table \ref{tab:bridge_SBC_PR_RBC_RB_sizes} gathers the sizes of the different reduced bases and the computation time to run the offline stage needed by Level 1 reduction. All simulations considered in this work were run on a 4-core laptop (with a 3.5 GHz Intel CPU and 16 GB RAM). We also provide the PR-RBC offline cost, though this cost is amortized over the many Level 1 - Level 2 queries.

\begin{table}[htb]
\begin{center}\begin{tabular}{|p{10cm}|p{5cm}|} \hline
	Size of training set $\Xi_{\rm o-t}$ & $41$ \\ \hline
	Size of port spaces & \begin{tabular}{@{}l@{}}
                   			$10$ for reference ports  $1$ and $3$ \\
                   			$12$ for reference port  $2$\\
                				\end{tabular} \\ \hline 
        Size of bubble spaces for port mode lifting & $6$ \\ \hline 
       Size of bubble space for inhomogeneity for archetype component $4$ & $10$  \\ \hline
       Computation time to run PR-RBC offline stage & $24.3\; s$ \\ \hline 
\end{tabular}\end{center}
\caption{PR-RBC reduced bases sizes for elastodynamics bridge}
\label{tab:bridge_SBC_PR_RBC_RB_sizes}
\end{table}%

\begin{figure}[htb]
	\begin{center}
	\includegraphics[scale=0.33]{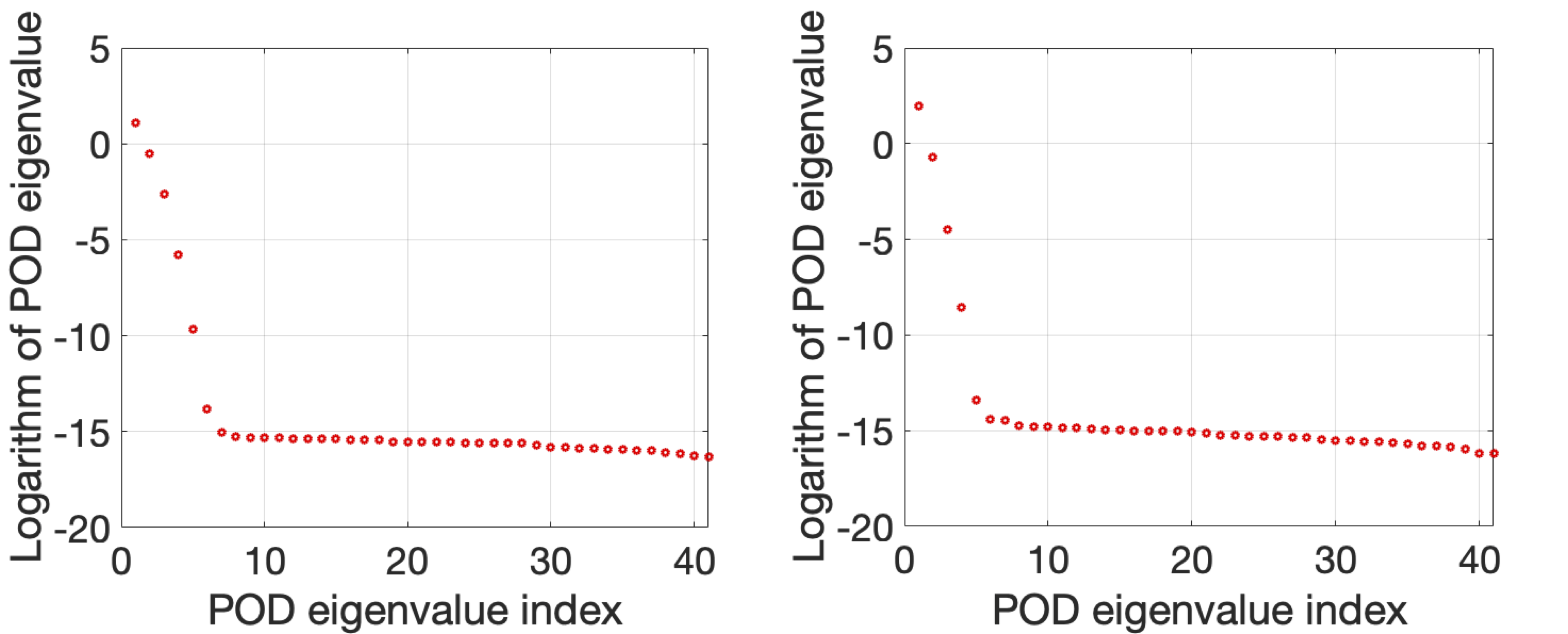}
	\caption{Logarithm of the POD eigenvalues for constructing the reduced space for port mode liftings for the 1st port mode retained for the reference port number $2$ and (left plot) archetype component $3$, and (right plot) archetype component $4$}
	\label{fig:bridge_port_4_port_mode_lifting_POD_eigenvalue_error}	
	\end{center}	
\end{figure}

\begin{figure}[htb]
	\begin{center}
	\includegraphics[scale=0.3]{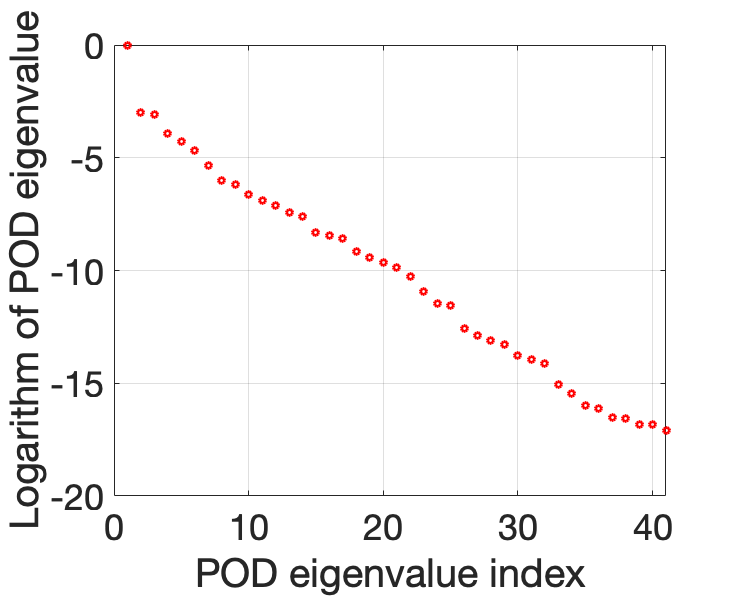}
	\caption{Logarithm of the POD eigenvalues for constructing the bubble space for inhomogeneity for archetype component number $4$}
	\label{fig:bridge_comp_4_POD_eigenvalue_error}	
	\end{center}	
\end{figure}

As an example, figure \ref{fig:bridge_port_4_port_mode_lifting_POD_eigenvalue_error} shows the decrease of the POD eigenvalues for the construction of reduced space for port mode liftings for the first port mode retained for the reference port number $2$ within the two archetype components number $3$ and $4$. We chose the corresponding reduced bubble spaces for port mode liftings to be of size equal to $6$. In figure \ref{fig:bridge_comp_4_POD_eigenvalue_error}, we give the convergence of the POD eigenvalues for construction of bubble space for inhomogeneity for archetype component number $4$, which is then chosen to be of size equal to $10$.





\subsection{Global System}

\begin{figure}[htb]
	\begin{center}
	\includegraphics[scale=0.15]{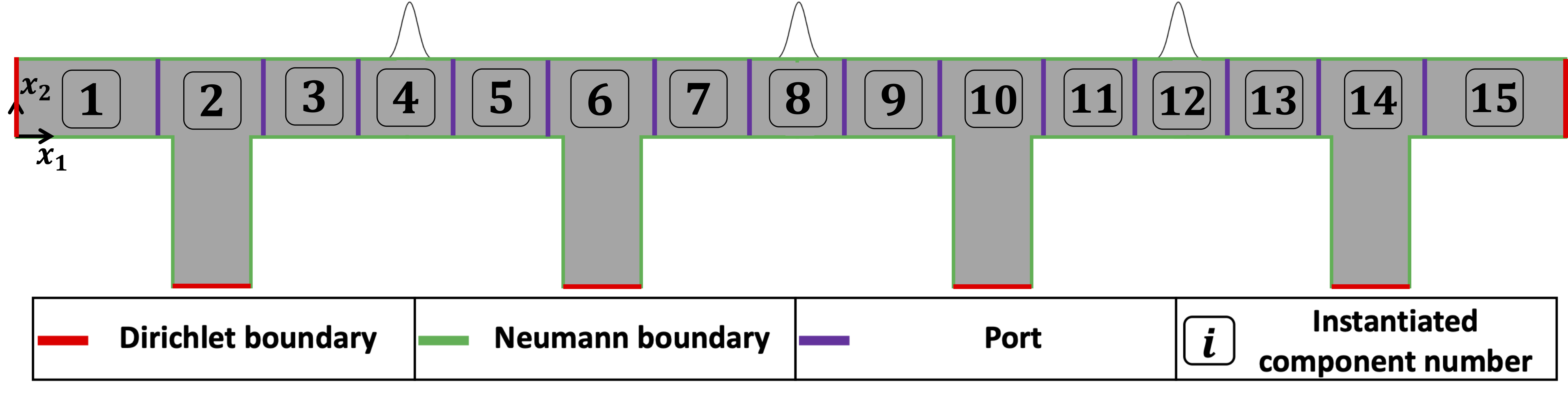}
	\caption{Not-to-scale representation of the global system for elastodynamics bridge: the Gaussian curves on the top boundary indicate the potential locations of the load}
	\label{fig:bridge_PR_RBC_glob_model_time}	
	\end{center}	
\end{figure}

We consider the global system consisting of $15$ instantiated components presented in figure \ref{fig:bridge_PR_RBC_glob_model_time}. The mapping of each instantiated component of the global system to the corresponding archetype is given in table \ref{tab:bridge_mapping_comp}. The time-domain variational problem for the global system has a total of $63$ parameters: the Young's modulus and the two Rayleigh damping coefficients for each component, in addition to $6$ parameters for each of the three beams where we can have a load ($5$ of these $6$ parameters appear naturally in equation (\ref{eq:nh_NBC}), and the last signals the existence or not of the load). As mentioned above, the frequency-domain variational problem has either $4$ or $8$ parameters for the archetype components, and $7$ or $11$ parameters for the archetype bi-components defining the reference ports. Those sizes are lower than the size of the parameter space of the global domain (equal to $63$) and thus it shows how the two-level PR-RBC method reduces the effective dimensionality of the parameter spaces considered in the variational problems. This reduction of the size of the parameter spaces is even more enhanced as larger global domains with more instantiated components are considered, since the sizes of the parameter spaces considered for the archetype components and reference ports do not change, while the size of the global parameter space will increase.

\begin{table}[htb]
\begin{center}\begin{tabular}{|p{7cm}|p{7cm}|} \hline
	Instantiated component number in global system & Archetype component number \\ \hline
	$1,15$ & $1$ \\ \hline 
        $2,6,10,14$ & $2$ \\ \hline
        $3,5,7,9,11,13$ & $3$ \\ \hline
        $4,8,12$ & $4$ \\ \hline
\end{tabular}\end{center}
\caption{Instantiated components to archetypes mapping}
\label{tab:bridge_mapping_comp}
\end{table}%



For the time-domain problem approximation, the simulation time $T_{\rm final}$ is taken equal to $800\; T_{\rm ref}$ such that we allow enough time for wave reflection. The size of the reduced space $X_{\rm RB}$ constructed by Strong Greedy approach (Algorithm \ref{alg:SG}) is fixed by imposing $\epsilon=10^{-5}$ as a threshold for the relative error $e_i$, $1\leq i\leq N$. The number of time steps $N_t$, and equivalently the step-size $\Delta t$, are fixed based on the convergence of the quantity:
\begin{equation}
\delta_{\Delta t}\equiv\frac{\max\limits_{1\leq j\leq N_t/2}\Big|\Big|\urbfd^{2\;j}(\mu)-\urbfderr^j(\mu)\Big|\Big|_{H^1(\Omega)}} {\max\limits_{1\leq j\leq N_t} \Big|\Big|\urbfd^j(\mu)\Big|\Big|_{H^1(\Omega)}} \ ,
\end{equation}
\noindent as $N_t$ increases. As expected for the second order Newmark-$\beta$ scheme ($\beta_t=\frac{1}{4}$, $\gamma_t=\frac{1}{2}$), the order of convergence in time-discretization is equal to $p=2$. Hence, a normalized Richardson's extrapolation-based error indicator is given by:
\begin{equation}
\epsilon_{\Delta t}\equiv\frac{1}{\max\limits_{1\leq j\leq N_t} \Big|\Big|\urbfd^j(\mu)\Big|\Big|_{H^1(\Omega)}}\;\frac{\max\limits_{1\leq j\leq N_t/2}\Big|\Big|\urbfd^{2\;j}(\mu)-\urbfderr^j(\mu)\Big|\Big|_{H^1(\Omega)}} {2^p-1} \ .
\end{equation}
\noindent The number of time steps $N_t$ is fixed such that we impose $\epsilon_{\Delta t}\leq10^{-3}$. 

For purposes of presentation, we consider a global parameter $\mu_{\rm example}$ such that the Young's modulus of all the components is taken equal to the nominal value, the Rayleigh damping coefficients are equal to half of their maximum value, and three loads are applied on the structure such they are centered at the middle of components number $4$, $8$ and $12$. The load parameters for the global parameter $\mu_{\rm example}$ are detailed in table \ref{tab:bridge_load_param}. We also consider $10$ randomly sampled parameters $\Xi_{\rm o}\equiv\{\mu_{{\rm rand} \, i}, 1\leq i\leq 10\}$. For each of these sampled parameters, we allow any combination regarding the existence or not of the load on components $4$, $8$ and $12$ (in total there are $2^3-1=7$ possible combinations for the load existence), and the temporal parameter $\sigma_t$ for the load applied on each of components $4$, $8$ and $12$ is sampled uniformly from the interval $[0.75\;\sigma_t^{\rm ref},1.25\;\sigma_t^{\rm ref}]$.

\begin{table}[htb]
\begin{center}\begin{tabular}{|p{4cm}|p{3cm}|p{3cm}|p{3cm}|} \hline
	Load on component \# & $4$ & $8$ & $12$ \\ \hline
	$\sigma_x$ & $0.02\; m$ & $0.03\; m$ & $0.04\; m$ \\ \hline 
       $F$ & $-20\;\frac{\bar{E}}{T_{\rm ref}}$ & $-15\;\frac{\bar{E}}{T_{\rm ref}}$ & $-10\;\frac{\bar{E}}{T_{\rm ref}}$ \\ \hline 
        $\sigma_t$ & $0.75\;\sigma_t^{\rm ref}(=12\; T_{\rm ref})$ & $\sigma_t^{\rm ref}(=16\; T_{\rm ref})$ & $1.25\;\sigma_t^{\rm ref}(=20\; T_{\rm ref})$ \\ \hline 
        $c_{\rm friction}$ & $0.7$ & $0.6$ & $0.5$ \\ \hline
\end{tabular}\end{center}
\caption{Load parameters for the global parameter $\mu_{\rm example}$}
\label{tab:bridge_load_param}
\end{table}%

\begin{figure}[htb]
	\begin{center}
	\includegraphics[scale=0.34]{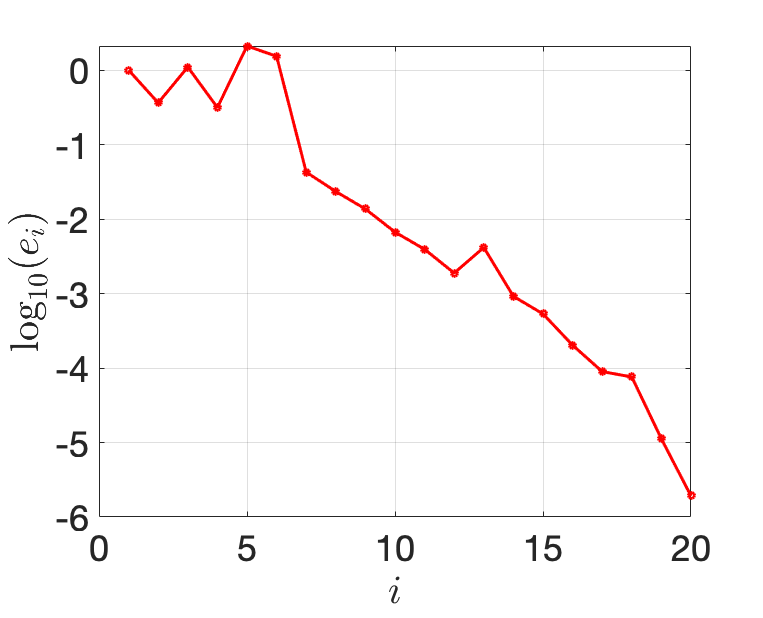}
	\caption{Strong Greedy algorithm (Algorithm \ref{alg:SG}) convergence for the elastodynamics bridge and the global parameter $\mu_{\rm example}$; here $n_\omega=41$}
	\label{fig:greedy_online_bridge_PR_RBC}	
	\end{center}	
\end{figure}

\begin{figure}[htb]
	\begin{center}
	\includegraphics[scale=0.35]{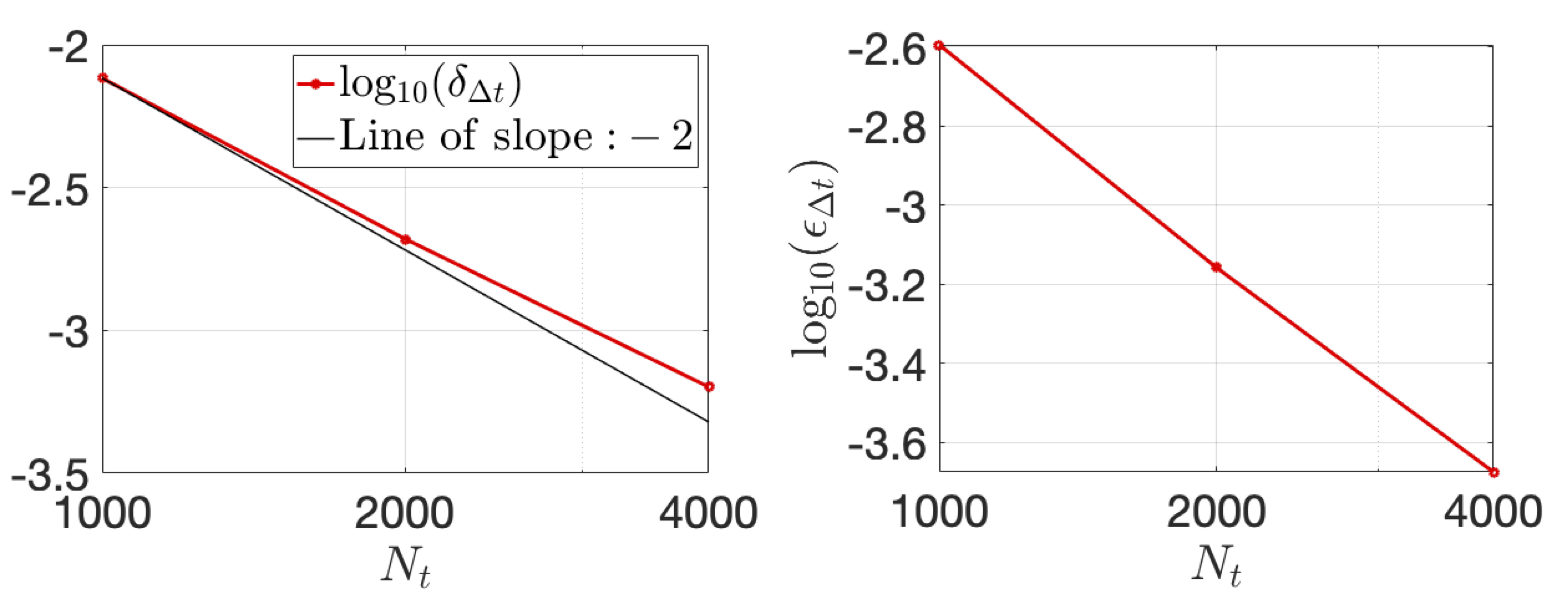}
	\caption{Convergence of $\delta_{\Delta t}$ and $\epsilon_{\Delta t}$ with $N_t$ for the global parameter $\mu_{\rm example}$}
	\label{fig:cv_time_2e3}
	\end{center}	
\end{figure}

The size of $X_{\rm RB}$ obtained by imposing $\epsilon=10^{-5}$ in algorithm \ref{alg:SG} for the global parameter $\mu_{\rm example}$ is $N=20$ as shown in figure \ref{fig:greedy_online_bridge_PR_RBC} giving the convergence of the Strong Greedy algorithm in terms of the relative error $e_i$, $1\leq i\leq N$. The second order convergence of $\delta_{\Delta t}$ with the number of time steps $N_t$ is verified for $\mu_{\rm example}$ as shown in the left plot of figure \ref{fig:cv_time_2e3}. Based on the right plot of figure \ref{fig:cv_time_2e3}, in which we present the convergence of $\epsilon_{\Delta t}$ with $N_t$, the number of time steps is taken as $N_t=2\times10^3$. We also verified the second order convergence of $\delta_{\Delta t}$ with $N_t$ for the $10$ randomly sampled parameters in $\Xi_{\rm o}$, and for $9$ of these cases the criterion $\epsilon_{\Delta t}\leq 10^{-3}$ requires $N_t=2\times10^3$, while for the remaining case, it requires $N_t=10^3$. Note that using explicit finite difference schemes generally results in a larger number of required time steps for stability reasons. For instance, using the explicit central difference scheme for the global parameter $\mu_{\rm example}$ requires at least $N_t=5\times10^5$ time steps for stability (corresponding to $250\times (2\times10^3)$, with  $N_t=2\times10^3$ being the number of time steps needed when using the implicit mid-point rule).

\begin{figure}[htb]
	\begin{center}
	\includegraphics[scale=0.45]{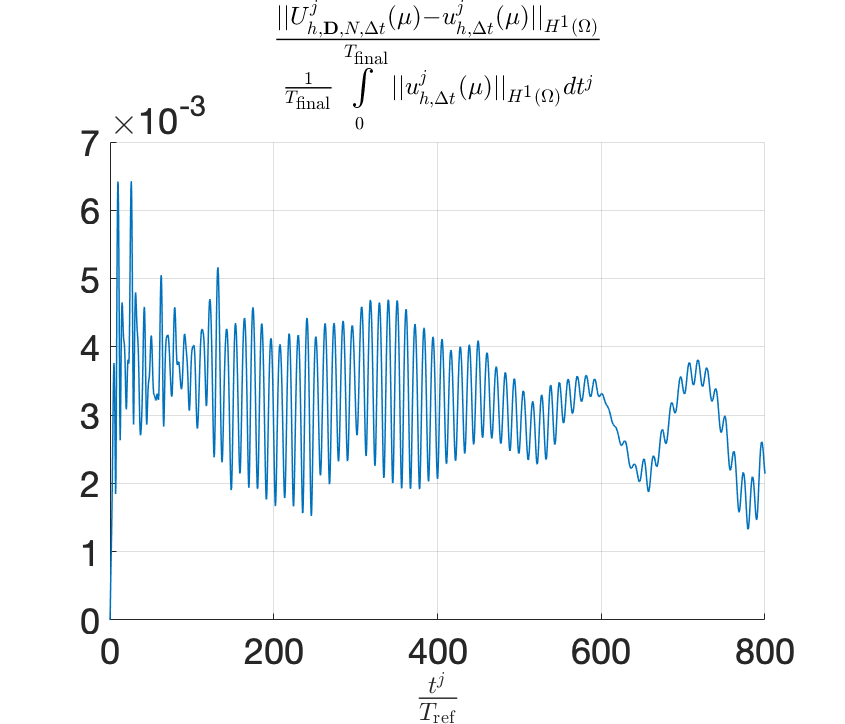}
	\caption{Time-domain relative error for the elastodynamics bridge and the global parameter $\mu_{\rm example}$}
	\label{fig:err_FE_PR_RBC_by_avg_time_bridge_PR_RBC}	
	\end{center}	
\end{figure}

Table \ref{tab:bridge_SBC_RB_online} gathers the computation time to estimate $\urbfd^j(\mu)$ , $1\leq j\leq N_t$, averaged over $\mu_{\rm example}$ and the $10$ randomly sampled parameters $\Xi_{\rm o}$. The size of the full $\PP_2$ FE approximation space is $\calN^h=11376$ and computing one full FE simulation with $N_t=2\times 10^3$ takes $1.77\; min$ on average (without performing the Richardson extrapolation). For confirmation purposes, we consider the time-domain relative error between the two-level PR-RBC solution and the FE approximation for a global parameter $\mu$,
\begin{equation*}
\frac{||\urbfd^j(\mu)-\ufefd^j(\mu)||_{H^1(\Omega)}}{\frac{1}{T_{\rm final}}\;\int\limits_0^{T_{\rm final}} ||u_{h,\Delta t}^j(\mu)||_{H^1(\Omega)} dt^j} \ , 1\leq j\leq N_t\ .
\end{equation*}
\noindent Figure \ref{fig:err_FE_PR_RBC_by_avg_time_bridge_PR_RBC} gives the evolution of the relative error for the global parameter $\mu_{\rm example}$ and we can verify that it is well below $1\%$, confirming the sufficiently refined time discretization thanks to the criterion $\epsilon_{\Delta t}<10^{-3}$, and the sufficiently rich reduced space $X_{\rm RB}$ thanks to the strong greedy criterion imposed with $\epsilon=10^{-5}$  in algorithm \ref{alg:SG}. Considering other relative errors, such as a normalization by the maximum value of $||u_{h,\Delta t}^j(\mu)||_{H^1(\Omega)}$ for $1\leq j\leq N_t$, gives very similar results.  Hence, having $n_\omega\leq50$ is enough to obtain sufficiently accurate approximations. By consequence, we typically have $N_t\gg n_\omega$ as mentioned in section \ref{subsec:op_count}. Note that the observed waves correspond well to dispersive flexural waves \cite{BhouriThesis}.

\begin{table}[htb]
\begin{center}\begin{tabular}{|p{12cm}|p{1cm}|} \hline
	PR-RBC online stage called $n_\omega=41$ times & $3.32 \;s$ \\ \hline
	Strong greedy (Algorithm \ref{alg:SG}) & $0.22\; s$  \\ \hline 
        Time marching (for $N_t=500,1000,2000,4000$ to perform the Richardson extrapolation) & $0.16\; s$   \\ \hline 
        Total computation time to estimate $\urbfd^j(\mu)$ , $1\leq j\leq N_t$ & $3.70\; s$ \\ \hline 
\end{tabular}\end{center}
\caption{Computation time of the proposed two-level reduced basis method for elastodynamics bridge}
\label{tab:bridge_SBC_RB_online}
\end{table}%

In conclusion, the two-level reduction approach has a computation cost $28.7$ times lower than the FE simulation for this example. For instance, in the context of Simulation Based Classification for Structural Health Monitoring, construction of datasets of size of the order of $10^4$ are often required to obtain satisfactory classification results \cite{BhouriThesis}. Conducting such a task for this bridge example using the two-level PR-RBC approach has a total computation time of $10.28$ hours (taking into account the computational cost of the PR-RBC offline stage, which is run only once), as opposed to an estimated $12.3$ days using the full FE approximation. 
 

\textbf{Acknowledgement: }
This work was supported by the ONR Grant [N00014-17-1-2077] and by the ARO Grant [W911NF1910098]. We would like to thank Dr. Tommaso Taddei and Professor Masayuki Yano for the helpful software they provided us with.






\newpage
\bibliographystyle{elsarticle-num}
\section*{\refname}
\bibliography{2_step_PR_RBC_final.bib}





\end{document}